\documentclass[leqno,11pt]{amsart}
\usepackage{latexsym,amsfonts,amsmath}
\usepackage{pb-diagram}

\newtheorem{theorem}{Theorem}[section]
\newtheorem{cor}[theorem]{Corollary}

\newtheorem{lemma}[theorem]{Lemma}
\newtheorem{prop}[theorem]{Proposition}
\theoremstyle{definition}
\newtheorem{definition}{Definition}[section]
\newtheorem{example}{Example}
\theoremstyle{remark}
\newtheorem{remark}{Remark}[section]
\numberwithin{equation}{section}
\newcommand{\nc}{\newcommand}
\newcommand{\C}{{\mathbb C}} \newcommand{\R}{{\mathbb R}}
\newcommand{\Z}{{\mathbb Z}}
 \newcommand{\Q}{{\mathbb Q}}
\newcommand{\CA}{\mathcal{A}}
\newcommand{\CB}{\mathcal{B}} 
 
\newcommand{\vol}{\operatorname{vol}}

\renewcommand{\c}{\mathfrak{c}}
\renewcommand{\a}{\mathfrak{a}}

\newcommand{\res}{\operatornamewithlimits{Res}}

\newcommand{\bff}{\mathbf{f}(z)}

\nc{\at}{\C[\ga]}
\nc{\rat}{\C_\A[\ga]}
\nc{\ga}{\mathfrak{a}}
\nc{\gac}{\ga_\C}
\nc{\gd}{\mathfrak{g}}

\nc{\s}{\mathfrak{s}}
\nc{\gv}{\mathfrak{v}}
\nc{\gt}{\mathfrak{t}}
\nc{\dga}{\mathfrak{a}^*}
\nc{\dgd}{\mathfrak{g}^*}
\nc{\dgt}{\mathfrak{t}^*}
\nc{\gma}{\Gamma_\ga}
\nc{\gmd}{\Gamma_\gd}
\nc{\gmt}{\Gamma_\gt}
\nc{\dgma}{\Gamma_\ga^*}
\nc{\dgmd}{\Gamma_\gd^*}
\nc{\dgmt}{\Gamma_\gt^*}
\nc{\scd}{{\sc scd}}
\nc{\gc}{\mathfrak{c}}
\nc{\A}{\mathfrak{A}}
\nc{\B}{\mathfrak{B}}
\nc{\spt}{\mathrm{Sec}}
\nc{\bA}{\mathrm{BInd}(\A)}
\nc{\bAF}{\mathrm{BInd}(\A, F)}
\nc{\bB}{\mathrm{BInd}(\B)}
\nc{\con}{\mathrm{Cone}}
\nc{\consing}{\mathrm{Cone}_\mathrm{sing}}
\nc{\gcp}{{\bar\gc}^\perp}
\nc{\comp}{U(\A)}
\nc{\tensor}{\otimes}
\nc{\ra}{\rightarrow}
\nc{\isom}{\cong}
\nc{\dvol}{d\mathrm{vol}}
\nc{\sumn}{\sum_{i=1}^n}
\nc{\sumr}{\sum_{i=1}^r}
\nc{\lr}[2]{\langle#1,#2\rangle}
\nc{\dt}[2]{#1=1,\dots,#2}
\nc{\fs}[1]{F^\sigma_{(#1)}}
\nc{\epu}{\epsilon_{\mathrm{u}}}
\nc{\epl}{\epsilon_{\mathrm{l}}}
\nc{\bss}{\mathrm{Basis}}
\nc{\permr}{\mathcal{S}_r}
\nc{\sign}{\mathrm{sign}}
\nc{\mpl}{\mathrm{MP}_\lambda}
\nc{\flaga}{\mathcal{FL}(\A)}
\nc{\flagat}{\mathcal{FL}(\A)_\tau}
\nc{\flagaut}{\mathcal{FL}(\Au)_\tau}
\nc{\flagau}{\mathcal{FL}(\Au)}
\nc{\flag}{\mathcal{FL}}

\nc{\pone}{\mathbb{P}^1}
\nc{\cone}{\mathrm{cone}}
\nc{\ar}{\rightarrow}
\nc{\bl}{{\boldsymbol{\lambda}}}
\nc{\bgam}{\boldsymbol{\gamma}}
\nc{\bn}{\mathbf{n}}
\nc{\solp}{\mathrm{sol}_{\A,\bl}(\bn)}
\nc{\bfb}{\mathbf{b}}
\nc{\sigm}{\bar\sigma}
\nc{\hc}{H_\mathrm{c}}
\nc{\hbm}{H^{\mathrm{BM}}}
\nc{\prbl}{p^\R_\bl}
\nc{\gC}{\mathfrak{C}}
\nc{\vb}{t}
\nc{\ve}{\boldsymbol{\epsilon}}
\nc{\bq}{\mathbf{q}}
\newcommand{\Cone}{\mathrm{Cone}}

\nc{\Eq}{\mathrm{Part}}

\nc{\chuze}{\mathrm{Choose}}
\nc{\FF}{\mathcal{FL}(\A,\xi)}

\nc{\FFF}{\mathcal{FL^+}(\A,\xi)}
\nc{\iml}{\mathrm{im}(L)}
\nc{\trop}{t}
\nc{\tS}{\widetilde{\mathrm{Eq}}}
\nc{\sol}{\mathrm{sol}}
\nc{\ts}{\widetilde{\mathrm{sol}}}
\nc{\vkf}{\boldsymbol{\kappa}^F}
\nc{\tD}{\widetilde{D}}
\nc{\vsigma}{{\boldsymbol{\sigma}}}
\nc{\const}{\mathrm{const}}
\nc{\cp}{{\mathrm{comp}}}
\nc{\ZZ}{\widehat{Z}}
\nc{\ZZF}{\widehat{Z}^F}
\nc{\JK}{\mathrm{JK}}
\nc{\overl}{\overset{\bl}<}

\newcommand{\rata}{\C_\A[\ga]}

\newcommand{\gr}{W}
\renewcommand{\gg}{\mathfrak{g}}
\newcommand{\ggy}{\mathfrak{g}}
\newcommand{\Crit}{\operatorname{Crit}}
\newcommand{\Torus}{\operatorname{T}_\C}
\newcommand{\Toric}{\operatorname{Tvar}_g}
\newcommand{\Hull}{\Pi}
\newcommand{\Spec}{\operatorname{Spec }}
\newcommand{\Hesszt}{\operatorname{Hess}_{z,\gt^*}}
\newcommand{\Hesszh}{\operatorname{Hess}_{z,\gt^*}}

\newcommand{\tres}{\mathrm{TorRes}}
\newcommand{\sing}{\mathrm{sing}}
\newcommand{\lra}{\longrightarrow}
\newcommand{\ggg}{\mathfrak{g}}
\newcommand{\gaz}{\ga_\Z}
\newcommand{\gtz}{\gt_\Z}
\newcommand{\ggz}{\ggg_\Z}

\newcommand{\tB}{\mathfrak{M}}
\newcommand{\M}{\mathfrak{M}}
\newcommand{\gamm}{\gamma}
\newcommand{\gwv}{V}
\newcommand{\gwt}{\gt^*}

\renewcommand{\ss}{l_\gamma}
\newcommand{\Map}{\pi}
\newcommand{\Tap}{\iota}
\newcommand{\tilCrit}{\widetilde{\Crit}}
\newcommand{\ru}{\underline{r}}
\newcommand{\nun}{\underline{n}}
\newcommand{\du}{\underline{d}}
\newcommand{\gau}{\underline{\ga}}
\newcommand{\ggu}{\underline{\gg}}
\newcommand{\gtu}{\underline{\gt}}
\newcommand{\Bu}{\underline{\B}}
\newcommand{\betau}{\underline{\beta}}
\newcommand{\alphau}{\underline{\alpha}}
\newcommand{\Au}{\underline{\A}}
\newcommand{\kappau}{\underline{\kappa}}
\newcommand{\deltau}{t}
\newcommand{\xiu}{\underline{\xi}}
\newcommand{\Fu}{\underline{F}}

\newcommand{\piu}{\underline{\pi}}
\newcommand{\iotau}{\underline{\iota}}
\newcommand{\gcu}{\gc}
\renewcommand{\gcd}{\bar{\gc}^\perp}
\newcommand{\idef}{{{\mathcal I}(z)}}

\newcommand{\bftz}{{\mathbf{f}(\tz)}}
\newcommand{\Face}{\Phi}
\newcommand{\TT}{\widehat T(\xi)}
\newcommand{\aaa}{\gt^*}
\newcommand{\Aff}{\mathrm{Aff}}
\newcommand{\II}{\mathrm{Im}'}
\newcommand{\tra}{\ga^*_{[\tau]}}
\newcommand{\xig}{\psi}
\newcommand{\xia}{\Psi}
\newcommand{\CH}{\Delta}
\newcommand{\PPu}{\underline{\Pi}}
\newcommand{\PPo}{\overline{\Pi}}
\newcommand{\tz}{{\tilde z}}
\newcommand{\themu}{\rho}

\nc{\bR}{\mathrm{RInd{\B}}}

\title{Mixed toric residues and tropical degenerations}
\author{Andr{\'a}s Szenes and Mich{\`e}le Vergne} 
\begin{document}
\maketitle

\setcounter{section}{-1}
\section{Introduction} \label{sec:intro}
This paper is a follow-up to our paper \cite{tor1}, where we prove a
conjecture of Batyrev and Materov, the Toric Residue Mirror Conjecture
(TRMC). Here we extend our results, and show that they imply a
generalization of this conjecture, the Mixed Toric Residue Mirror
Conjecture (MTRMC), which is also due to Batyrev and Materov
\cite{BMmixed}.

Roughly, these conjectures state that the generating function of
certain intersection numbers of a sequence of toric varieties
converges to a rational function, which can be obtained as a
finite residue sum on a single toric variety.  We first recall the
TRMC in some detail. We start with an integral convex polytope
$\Pi^\B$ in a $d$-dimensional real vector space $\gt$ endowed with
a lattice of full rank $\gtz$; we assume that the polytope
contains the origin in its interior. Let the sequence
$\B=[\beta_1,\beta_2,\ldots,\beta_n]$ be the set of vertices of
this polytope, ordered in an arbitrary fashion.  One can associate
a $d$-dimensional polarized toric variety $(V^\B, L^\B)$ to this
data in the standard fashion \cite{tor1}.

There is another way to obtain toric varieties from this data, which
generalizes the mirror duality of polytopes introduced by Batyrev
\cite{Bat1}.
Consider the sequence $\A=[\alpha_1,\alpha_2,\ldots,\alpha_n]$,
which is the Gale dual of $\B$ (cf. \S\ref{sec:gale} for the
construction). This is a sequence of integral vectors in the dual
$\ga^*$ of a certain $r=n-d$-dimensional vector space $\ga$, which
is also endowed with a lattice of full rank: $\gaz$; in this setup
the sequence $\A$ spans a strictly convex cone $\Cone(\A)$.  The
simplicial cones generated by $\A$ divide  $\Cone(\A)$ into open
chambers. Each chamber corresponds to a $d$-dimensional orbifold
toric variety $V_\A(\gc)$ (cf. \cite{Danilov}). An integral
element $\alpha\in\ga^*$ specifies an orbi-line-bundle $L_\alpha$
over this variety; denote the first chern class of $L_\alpha$ by
$\chi(\alpha)\in H^2(V_\A(\gc),\Q)$. For the purposes of this
introduction we assume that this correspondence induces the linear
isomorphisms
\[ \ga^*\isom H^2(V_\A(\gc),\R)\quad\text{and}\quad
    \ga\isom H_2(V_\A(\gc),\R)
\]

Now pick a chamber $\gc$ which contains the vector $\kappa=\sumn
\alpha_i$ in its closure: $\kappa\in\bar{\gc}$. To each element
$\lambda\in\gaz$, one can associate a moduli space $MP_\lambda$, the
so-called Morrison-Plesser space (cf. \cite{MP}), which is a
compactification of the space of those maps from the projective line
to the variety $V_\A(\gc)$ under which the image of the fundamental
class is $\lambda$:
\[ \{m:\pone\ar V_\A(\gc);\;m_*([\pone])=\lambda\}.\]
The varieties $\mpl$ are toric, and such that, again, to each integral
element $\alpha \in \ga^*$ one can associate a line bundle
$L_{\alpha}$ on $\mpl$; again we denote the corresponding chern class
in $H^2(\mpl)$ by $\chi(\alpha)$.  The space $\mpl$ is defined to be
empty, unless $\lr\alpha\lambda\geq0$ for every
$\alpha\in\gc$. The set of vectors satisfying this condition forms a
cone in $\ga$, which we denote by $\gcp$; this cone is called the {\em
polar} cone of $\gc$.

The construction also provides a Poincare dual class 
\[K_\lambda\in H^{2(\dim \mpl-d)}(\mpl,\Q)
\]
to the subspace of $\mpl$ of those maps which land in a generic
zero-section $Y$ of the line bundle $L_{\kappa}$.  When $V_\A(\gc)$ is
smooth, then $Y$ is a Calabi-Yau manifold.

To probe the class $K_\lambda$, we fix a homogeneous polynomial
$P(x_1,\dots,x_n)$ of degree $d$ in $n$ variables, and consider
the intersection numbers
\[\int_{\mpl}P(\chi(\alpha_1),\dots,\chi(\alpha_n))K_\lambda,\]
which are to be interpreted as analogs of numbers of rational
curves in $Y$ subject to certain conditions specified by the
polynomial $P$.

Now let $z_1,\dots,z_n\in\C^*$, and form the Laurent series
\begin{equation}
\label{formalseries} \sum_{\lambda\in\ga_\Z}
\int_{\mpl}P(\chi(\alpha_1),\dots,\chi(\alpha_n))K_{\lambda}
\prod_{i=1}^n z_i^{\lr{\alpha_i}\lambda}.
\end{equation}

The statement of the conjecture is that this series converges in a
certain domain of $(\C^*)^n$ depending on the chamber $\gc$, and the
limit is a rational function of the parameters $z_1,\dots,z_n$, which
is given by a {\em toric residue} on the polarized toric variety
$(V^\B,L^\B)$.  This limiting function does not depend on which
chamber $\gc$ with $\kappa\in\bar{\gc}$ we picked.

Toric residues for polarized toric varieties were introduced by Cox
\cite{Cox}.  Under appropriate conditions, a polarized toric variety
$(V,L)$ of dimension $d$, and a generic $(d+1)$-tuple of homogeneous
functions $\bff=(f_0,f_1,\dots,f_d)$ in $H^0(L,V)$ defines a residue
functional on the homogeneous ring 
\[ Q\longmapsto\tres_{\bff} Q. \] This functional vanishes on the
ideal generated by $\bff$, and depends rationally on $\bff$ and $Q$.
In our case, by definition, the space of sections $H^0(L^\B,V^\B)$ of
the polarizing line bundle $L^\B$ is spanned by
$\{e_{\beta};\,\beta\in\Pi^\B\cap\gtz\}$, where $e_\beta$ is a
character of the torus $\gt^*/\gtz^*$. In particular, the section
$f_0=1-\sumn z_ie_{\beta_i}$ and the $d$ sections $\{\sumn\lr
a{\beta_i}z_ie_{\beta_i};\,a\in\mathbf{a}\}$, where $\mathbf{a}$ is an
integral unimodular basis of $\gt^*$, generate an appropriate ideal if
the parameters $z_1,\dots, z_n$ are sufficiently generic. For $Q$ we
take the element $P(z_1e_{\beta_1},\dots,z_ne_{\beta_n})$ of degree
$d$ of the homogeneous ring.  Then the TRMC may be more precisely
formulated as follows: the rational function of the parameters $z_i$
obtained as the toric residue
$\tres_{\bff}P(z_1e_{\beta_1},\dots,z_ne_{\beta_n})$ may be expanded
in a domain of $(\C^*)^n$ depending on the chamber $\gc$, and the
resulting Laurent series  is given by \eqref{formalseries}.

In \cite{tor1} we gave a proof of this conjecture.  The central result
to which we reduce the proof of TRMC is a statement about the topology
of toric varieties: Let $\comp$ be the complement of the complexified
hyperplane arrangement in $\ga\tensor_\R\C$ induced by the sequence
$\A$. Note that the map $\chi:\ga^*\rightarrow H^2(V_\A(\gc),\R)$,
defined above, may be extended multiplicatively to a map from the
space of polynomials on $\ga$ to the even cohomology ring of
$V_\A(\gc)$.  We show that a certain compact real algebraic cycle
$Z(\xi)\subset\comp$ of dimension $r$,  depending on a parameter
$\xi\in\gc$, computes the intersection ring of the toric variety
$V_\A(\gc)$ in the following sense. For a polynomial $Q$ of degree $d$
on $\ga$, we have
\[ \int_{V_\A(\gc)} \chi(Q) = \int_{Z(\xi)}
\frac{Q\,\dvol}{\prod_{i=1}^n\alpha_i},\] where $\dvol$ is the
complexification of the volume form on $\ga$ induced by $\gaz$.
The cycle $Z(\xi)$ is given by the following real algebraic equations
\[  Z(\xi)=\left\{u\in\comp; \;
  \prod_{i=1}^n|\alpha_i(u)|^{\langle\alpha_i,\lambda\rangle}= 
e^{\langle\xi,\lambda\rangle}\text{ for all }\lambda\in\gaz
\right\}.
\]
The proof of this result given in \cite{tor1} uses a degeneration
technique related to the so-called {\em tropical geometry}
\cite{Viro,Sturm}. We will demonstrate the argument in a simple
example (cf. Example \ref{ex1}).


The Mixed TRMC is a similar naive rational curve counting formula
for higher codimensional complete intersections in toric
varieties. Here we start with a partition $\{D_1,\dots,D_l\}$ of
the index set $\{1,\dots,n\}$ of the sequence $\A$, which
satisfies the following condition: there exists a chamber $\gc$ of
$\A$ such that the vectors $\theta_s=\sum_{i\in D_s}\alpha_i$,
$\dt sl$ are in the closure of $\gc$.  Again, a common smooth
generic intersection of zero-sections of the line bundles
$L_{\theta_s}$, $\dt sl$, is a Calabi-Yau submanifold in
$V_\A(\gc)$.

Batyrev and Materov \cite{BMmixed} formulate their conjecture, the
MTRMC, for so-called {\em nef-partitions}; this notion is a
reformulation of the condition we described above.

The structure of the conjectured formula in the mixed case is similar
to that of the TRMC. We form a series of the form of
\eqref{formalseries} with a suitably adjusted class $K_\lambda$ (see
\eqref{interjk} and \eqref{theseries}). Again we need to show that
this series converges to a rational function which may be expressed as
a toric residue. In this case, the toric residue corresponds the the
so-called Cayley polytope instead of the polytope $\Pi^\B$
(cf.\cite{BMmixed}). This is a polytope of dimension greater than $d$;
it is constructed from $\Pi^\B$ and the data of a ``nef''-partition.

The definitions of the series involved in the Batyrev-Materov
conjectures are given by expansions of rational functions, and are
relatively simple.  It is unclear how to extend our approach to the
computation of the Givental $J$-series \cite{Giv}.  The conjectures of
Batyrev-Materov do, however, provide another contribution to the
mirror symmetry phenomena.

In this paper, we present a proof of this conjecture. We show, that
the MTRMC can again be reduced to our integration
formula for toric varieties. Besides this reduction, the paper
contains several improvements over \cite{tor1}. These include
\begin{itemize}
\item a more invariant definition of the relevant maps in our
  construction, which does not rely on the choice of a special
  $\Z$-basis of $\gaz$ (cf. \eqref{diagram2})
\item a  more natural and less restrictive condition of regularity
  on the element $\xi$, which guarantees our statements on
  $Z(\xi)$; we call this condition {\em flag-regularity}
  (cf. Definitions \ref{regular}, \ref{defsplus}).
\item a more general and simpler Laurent expansion formula which
  is valid for any chosen basis of $\gaz$ (cf. Proposition \ref{sum}).
\end{itemize}

Now we review the contents of the paper. In \S\ref{sec:localform} we
recall the definition of the toric residue and describe a localized
formula for it. We apply the idea from \cite{tor1}, which
allows us to write the localized formula as a sum of the values of a
rational function on $\ga$ over the finite set of solutions of a
system of binomial equations. Each of the equations has the form
$\prod_{i=1}^n \alpha_i^{\nu_i}=t$, where $\nu_i\in\Z$ for $\dt in$,
and $t$ is a complex parameter.
In \S\ref{sec:summary} we review the main results of \cite{tor1} with
the improvements described above.
Finally, in \S\ref{sec:mixed}, we turn to the proof of the MTRMC of
Batyrev and Materov. We start by translating the Cayley construction
into the Gale dual language of the sequence $\A$. This means the
following: we add one dimension to the space $\ga$. We denote the
extra coordinate on $\ga$ by $t$. Next we replace $\A$ with the new
collection formed by the vectors $[\alpha_i,0],i=1,\ldots,n$ together
with $l$ new vectors $[-\theta_k,t],k=1,\ldots,l$.  Applying our
results from the previous section to this augmented sequence, we
quickly arrive at the proof of the conjecture.

We end the paper with an appendix, where we justify our use of the
local formula, Proposition \ref{tlocal}, for the toric residue in the
mixed case. Here this is rather important, since our toric residue is
taken in the Cayley toric variety corresponding to the Cayley
polytope, but the deformation parameters vary in a smaller dimensional
space. As a result, we cannot claim to be taking generic values of the
parameters. This makes checking the conditions of the Propositions
more difficult. It is necessary to justify why on the particular
subset the assumptions hold. In fact, we compute the exact domain of
the parameter space where the conditions hold.

We would like to note that, just as our previous result in
\cite{tor1}, we do not need reflexivity. Our formula holds in much
greater generality.

We were informed that a few months earlier  Kalle Karu \cite{karu} has
given a proof of this conjecture using different methods.

{\bf Acknowledgment}. We would like to thank Alicia Dickenstein for
her help.  We are grateful to the Erwin Schr\"odinger
Institute for excellent working conditions. The first author would
like to acknowledge the support of the Hungarian Science Foundation (OTKA).

\section{Local formula for the toric residue}
\label{sec:localform}

\subsection{Notation and Data}
\label{sec:notation}
Let $V$ be a finite-dimensional real vector space of dimension
$d+1$ and denote by $\gwv_\C$ its complexification. Assume that
$\gwv$ is endowed with an integral structure, i.e. with a lattice
$\gwv_\Z$ of full rank. Then denote by $\Torus(\gwv)$ the
complexified torus $\gwv_\C/\gwv_\Z$, and by $\gwv_\Z^*$ the set
of linear forms on $\gwv$ which assume integral values on
$\gwv_\Z$. Thus
\[ \gwv_\Z^*=\{a\in\gwv^*;\; \lr a\mu \in\Z\text{ for all }\mu\in\gwv_\Z\}. \]
For each vector $\mu\in\gwv_\Z$, denote by $e_{\mu}$ the function
$e^{2\pi i\mu}$ on $V^*_\C$, which may also be considered as a
function on $\Torus(\gwv^*)$.

Further, fix a primitive element $g\in V^*_\Z$ and a sequence
$\tB=[\mu_1,\mu_2,\ldots, \mu_n]$ of vectors of the lattice $V_{\Z}$,
which generate $V_\Z$ over $\Z$, and satisfy $\langle g,\mu_i\rangle
=1$, $i=1,\dots,n$. Thus $\tB$ is a sequence of vectors in an affine
subspace of $V$ of codimension 1. Let $C=C(\tB)\subset V$ be the acute
$d+1$-dimensional polyhedral cone generated by $\tB$:
$C=\sumn\R^+\mu_i$. The convex hull $\Hull(\tB)$ of the elements of
$\tB$ is a convex polytope of dimension $d$, which serves as the base
of the cone $C$. The faces of the polyhedral cone $C$ are polyhedral
sub-cones of $C$ which have dimensions running from $0$ to $d+1$.
Thus, under our convention, the cone $C$ itself is a face.

Introduce the algebra $S(C)$ with basis given by the functions
$e_{\mu}$, $\mu\in V_\Z\cap C$.  This algebra is graded:
\[S(C)=\oplus_{k=1}^{\infty} S^k(C)\text{ with }
 S^k(C) = \bigoplus \C e_\mu,\; \mu\in V_\Z\cap C,\;
\lr{g}\mu=k.\] This algebra has a natural ideal
$I(C)=\oplus_{k=1}^{\infty} I^k(C)$, the so-called {\em dualizing ideal},
which is generated by the elements $e_{\mu}$ with $\mu$ in the
interior of the cone $C$.

The finitely generated algebra $S(C)$ defines an affine variety
$\Aff(C)=\Spec(S(C))$ of dimension $d+1$.  A point in $\Aff(C)$ is a
character $x$ of $S(C)$, i.e. a ring homomorphism $x:S(C)\ar\C$; we
will write $e_\mu(x)$ instead of $x(e_\mu)$ and we will use the
notation $x_i$ for $e_{\mu_i}(x)$, $\dt in$. The variety $\Aff(C)$ is a
closed cone, and it contains $\Torus(V^*)$ as a Zariski open subset.
Indeed, there is a one-to-one correspondence between points of
$\Torus(V^*)$ and characters $x$ of $S(C)$ such that $x_i\neq 0$ for
all $i=1,2,\ldots,n$.

The graded algebra $S(C)$ defines a projective toric variety $\Toric
(C)$ of complex dimension $d$, together with an ample line bundle
$L\to \Toric(C)$. Every element $f\in S^k(C)$ defines a section of the
line bundle $L^k$; we will denote this section by the same symbol $f$.
Each such section defines an hypersurface in $\Toric(C)$ via its set
of zeros. Denote by $H$ the subspace of $V$ orthogonal to our grading
vector $g$, with lattice $H_\Z=H\cap V_\Z$.  Then $H^*$ is identified
with $V^*/\R g$ and the complexified torus
$\Torus(H^*)=\Torus(V^*)/\C^*$ is naturally embedded in $\Toric(C)$ as
a Zariski open subset.

We choose a primitive element $\gamm\in V_\Z$ in the interior of $C$,
and denote by $\ss$ its degree: $\ss=\langle g,\gamm\rangle$.  Denote
by $\gt$ the quotient $V/\R\gamm$ endowed with the lattice which is
the image of $V_\Z$ in this quotient. Then $\gwt$ is naturally
embedded in $V^*$ as the subspace orthogonal to $\gamm$. We also have
$\gt^*_{\Z}= V^*_{\Z}\cap\gt^*$.

Introduce the auxiliary vector space $\ggy=\oplus_{i=1}^n \R
\omega_i$ endowed with the lattice
$\ggy_\Z=\oplus_{i=1}^n \Z \omega_i$. Let $\pi:\ggy\ar\gt=V/\R
\gamm$ be the surjective map which sends $\omega_i$ to the image
of $\mu_i$ modulo $\R\gamm$. We introduce the notation $\beta_i$
for these image vectors $\Map(\mu_i)$; thus we have a sequence of
vectors
$${\B}=[{\beta_1},\ldots, {\beta_n}]\subset\gt.$$
The elements of this
sequence generate $\gt_\Z$ over $\Z$. Note that the natural inclusion
$\gt^*\hookrightarrow V^*$ induces the inclusion
$\Torus(\gt^*)\hookrightarrow \Torus(V^*)$. We make the obvious
observation that
\begin{equation}
  \label{obvious}
\text{if }w\in\gt^*,\text{ then }\langle \beta_i,w\rangle =\langle
\mu_i,w\rangle\text{ and }e_{\beta_i}(w)=e_{\mu_i}(w),\,\dt in.
  \end{equation}

Denote by $\ga$ the kernel of the map $\Map$, endowed as usual by the
integral structure inherited from $\ggy$; then we have the sequence
\begin{equation}
  \label{exact1}
0\to \ga\overset{\Tap}\to \ggy\overset{\Map}\to \gt\to 0,
\end{equation}
which is exact over $\Z$.

Considering the map $\gg\ar V$, which, similarly to $\pi$, associates
$\mu_i$ to $\omega_i$, and denoting its kernel by $W$, we can arrange
our real vector spaces in the following diagram:
\begin{equation}
\label{diagram}
\begin{diagram}
\node[3]{\ga}\arrow{s,l}{\iota}  \node{H}\arrow{s} \node{\R\gamm}\arrow{sw}
\\
\node{0}\arrow{e}\node{W}\arrow{e}\arrow{ne}\node{\gg}\arrow{e}\arrow{s,l}{\pi}
\node{V}\arrow{sw}\arrow{e}\arrow{s,r}{g}\node{0}
\\
\node[3]{\gt}\node{\R}
\end{diagram}
\end{equation}
Each of the vector spaces is endowed with an integral structure and
the exact sequences which appear on the diagram are exact over $\Z$ as
well. The vector sequences we have introduced so far are $\mu_i\in
V,\,\beta_i\in\gt,\, \omega_i\in\gg$, $\dt in$.

Note that if $\gamm=\sum_{i=1}^n r_i \mu_i$, we can write
$$\ga=\gr\oplus \R \sum_{i=1}^n r_i \omega_i.$$

\begin{lemma}\label{covering}
  Consider the map $H\ar\gt$ which is obtained by composing two arrows
  on the diagram \eqref{diagram}.  This map induces covering maps
  $$\Torus(H)\ar\Torus(\gt)\text{ and }\Torus(\gt^*)\ar\Torus(H^*)$$
  of order $\ss$.
\end{lemma}

\begin{proof}
This follows from the dual fact that the image of the vector $\gamma$
under the linear functional $g$ is $\ss$.

\end{proof}

\subsection{The toric residue}\label{sec:torres}
Let $z=(z_1,z_2,\ldots, z_n)$ be a generic vector in $\C^n$, and
introduce the element
$$f_z=\sum_{i=1}^n z_i e_{\mu_i}\in S^1(C).$$
Considered as a section
of the polarizing line bundle $L$ over $\Toric(C)$, it defines a
hypersurface $\{f_z=0\}\subset\Toric(C)$ varying with $z$.  For
each $a\in V^*$, we consider the derivative
$$f_{z,a}=\sum_{i=1}^n z_i \langle a,\mu_i\rangle e_{\mu_i}$$
of the function $f_z$ in the direction of the vector $a\in V^*$.
Again the functions $f_{z,a}$ are elements of $S^1(C)$.  Observe
that $f_{z,g}=f_z$.

Next, we choose a {\em $\Z$-basis} $ [a_1,\ldots,a_{d}]$ of
$\gt^*_\Z$; these induce $d+1$ sections $f_{z,0},f_{z,1},\ldots,
f_{z,d}$ of $L$: $f_{z,0}=f_z$ and $f_{z,j}=f_{z,a_j}$, $\dt jd$.

Now consider the following conditions:
\begin{eqnarray}\label{cond1}
& \{x\in\Toric(C);\; f_{z,j}(x)=0, \,j=0,\dots,d\}=\emptyset;\\
\label{cond2} &\{x\in\Toric(C);\; f_{z,j}(x)=0, \,j=1,\dots,d
\}\setminus\Torus(H^*)=\emptyset.
\end{eqnarray}
The first says the hyper-surfaces $\{f_{z,j}=0\}$, $j=0,1,\dots,d$, are in
generic position, and the second that the last $d$ intersect inside
the complexified torus of $\Toric(C)$. Note that the first condition
does not depend on our basis, while the second
depends only on the choice of the vector $\gamma$. Indeed, we could
write the two conditions as
\begin{eqnarray*}
& \{x\in\Toric(C);\; f_{z,a}(x)=0 \text{ for }a\in V^*\}=\emptyset;\\
&\{x\in\Toric(C);\; f_{z,a}(x)=0\text{ for }\lr a\gamma=0
\}\setminus\Torus(H^*)=\emptyset.
\end{eqnarray*}

Denote by $\idef$ the $S(C)$-ideal generated by
$\bff=\{f_{z,j},\,j=0,\dots,d\}$. According to \cite{CDS}, if condition
\eqref{cond1} is fulfilled, then the codimension of $\idef\cap
I^{d+1}(C)$ in $I^{d+1}(C)$ is exactly 1, and there is a
canonical functional defined by Cox \cite{Cox}:
\[ \tres_{\bff}: I^{d+1}(C)\longrightarrow \C\]
called the {\em toric residue}, which vanishes on $\idef\cap
I^{d+1}(C)$.

Now, following \cite[Proposition 2.6]{BMmixed}, we describe a local
formula for this toric residue.
For $z\in (\C^*)^n$ introduce the set
\[\Crit(\aaa,z)=\{x\in\Torus(H^*);\; f_{z,j}(x)=0, \,j=1,\dots,d \}.
\]
Note that if \eqref{cond2} holds, then this set coincides with
\[ \{x\in\Toric(C);\; f_{z,j}(x)=0, \,j=1,\dots,d \},\]
and that if \eqref{cond1} holds then $\Crit(\aaa,z)$ is finite.
We will need a variant of $\Crit(\aaa,z)$:
\[\tilCrit(\gt^*,z)=\{w\in\Torus(\gt^*);\; f_{z,j}(w)=0, \,j=1,\dots,d \},
\]
which is the inverse image of $\Crit(\gt^*,z)\subset \Torus(H^*)$ in
the covering group $\Torus(\gt^*)$. According to Lemma
\ref{covering}, this set is an $l_\gamma$-fold covering of
$\Crit(\gt^*,z)$.

Now consider the second derivatives of the function $f_z$:
\begin{equation}\label{secondder}
f_{z,jk}=\sum_{i=1}^{n} z_i \langle a_k,\mu_i\rangle \langle
a_j,\mu_i\rangle e_{\mu_i},\;j,k=1,\dots, d
\end{equation}
and let
\begin{equation}\label{defhess}
\Hesszh=\det(f_{z,jk})_{j,k\geq 1}
\end{equation}
be the determinant of the matrix formed by these second
derivatives. Then $\Hesszh$ is an element of $S^{d}(C)$, which depends
on $\aaa$ and $z$ only, and the product $f_z \Hesszh$ is an element of
$S^{d+1}(C)$.

Consider the following  non-degeneracy condition:
\begin{equation}\label{cond3}
 \Hesszh(x)\neq0\text{ for each }x\in\Crit(\aaa,z).
\end{equation}

Now we can formulate the localization theorem for toric residue:
\begin{prop}\label{tlocal}
  Let $Q\in I^{d+1}(C)$, and assume that $\gt^*$ and $z\in \C^n$ are
  such that conditions \eqref{cond1}, \eqref{cond2} and \eqref{cond3}
  hold.  Then $Q/(f_z\Hesszh)$ is a meromorphic function on
  $\Toric(C)$, whose poles avoid $\Crit(\aaa,z)$, and we have
\begin{equation}\label{localform}\tres_{\bff}Q=l_\gamma\sum_{w\in
\tilCrit(\aaa,z)} \frac{Q(w)}{f_z(w) \Hesszh(w)}.
\end{equation}
\end{prop}

This statement is almost identical to \cite[Proposition 2.6]{BMmixed},
and its proof quickly follows, for example, from \cite[Theorem
3.2]{BINOM}. The factor $l_\gamma$ appears for the following reason:
on the one hand, the basis $[g,a_1,\dots,a_d]$ is not a $\Z$-basis of
$V^*_\Z$, but a basis which spans a parallelepiped with volume
$l_\gamma$ in $V^*_\Z$; this changes the Hessian by a factor of
$(l_\gamma)^2$. On the other hand, we replaced the set $\Crit(\aaa,z)$ by
its $l_\gamma$-fold cover in the formula, and this gives us a factor of
$l_\gamma^{-1}$.

In Propositions \ref{denseopen} and \ref{possible} of the appendix, we
will describe a set $Z$ of elements $z\in \C^n$ such that conditions
\eqref{cond1}, \eqref{cond2} and \eqref{cond3} hold for $z\in Z$.

\subsection{Gale duality}
\label{sec:gale}
Now we turn to the description of the set $\Crit(\gt^*,z)$, which is a key
point of our work.  The situation is similar to that in our article
\cite{tor1}.

Consider the exact sequence \eqref{exact1}.  Denote by $\alpha_i\in
\ga^*$ the restriction of the coordinate function $\omega^i$ on $\gg$
to $\ga$. Let $\A$ be the sequence
$$\A:=[\alpha_1,\alpha_2,\ldots, \alpha_n].$$
Then the sequence $\A$
is defined to be {\em Gale dual} of the sequence $\B$.

Observe that the convex hull of the vectors in $\B$ contains the
origin in its interior, and as a result, the vectors of the Gale dual
sequence $\A$ lie in an open half-space of $\ga^*$. We called such a
sequence $\A$ {\em projective} in \cite{tor1}.

An element $u$ of $\ga$ is
written $u=\sum_{i=1}^n \alpha_i(u)\omega_i$ and for every $u\in \ga$,
we have
$$\sum_{i=1}^n \alpha_i(u)\mu_i\in \R \gamm.$$

The sequence $\A$ defines a hyperplane arrangement in $\ga_\C$ and
we denote by $U(\A)$ the complement of this hyperplane
arrangement:
$$U(\A)=\{u\in \ga_\C ;\, \alpha(u) \neq 0 \text{ for all }\alpha \in
\A\}.$$ For $\lambda\in \ga_\Z$, we introduce the rational functions
\begin{equation}
  \label{plam}
p_\lambda(u)=\prod_{i=1}^n \alpha_i(u)^{\langle
  \alpha_i,\lambda\rangle}\text{ and }
z^{\lambda}=\prod_{i=1}^n z_i^{\langle \alpha_i,\lambda\rangle}.
  \end{equation}
on $\ga_\C$ and $(\C^*)^n$, respectively.

Now, similarly to \cite{tor1}, we parameterize the set
$\tilCrit(\gt^*,z)\subset \Torus(\gt^*)$ by a finite subset $O(z,\A)$
in $\ga_\C$.

\begin{definition}\label{defo}
For $z\in(\C^*)^n$ define the set
 $$O(z,\A):=\{u\in U(\A);\; p_\lambda(u)=z^{\lambda}  \text{ for
all }\lambda \in \ga_\Z\}.$$
\end{definition}

\begin{lemma}\label{OZA}
  If $u\in O(z,\A)$, then there exists a unique element $w\in
  {\tilCrit}(\gt^*,z)$ such that
\begin{equation}\label{w} z_i e_{ \beta_i}(w)=\alpha_i(u)\;
\text{ for }\,\dt in.
\end{equation}
This correspondence between $O(z,\A)$ to $\tilCrit(\gt^*,z)$ is bijective.
\end{lemma}

\begin{proof}

  If $w$ is in ${\tilCrit}(\gt^*,z)$, then $f_{z,a}(w)=0$ for all $a\in
  \gwt$. Explicitly,
$$\sum_{i=1}^n z_i \langle a,\beta_i\rangle   e_{\beta_i}(w)=0
\text{ for all } a\in \gt^*.$$
 Thus
$$\sum_{i=1}^n z_i   e_{\beta_i}(w)\beta_i=0.$$

By definition of Gale duality, this condition means that there exists
$u\in \ga_\C$ such that $z_i e_{\beta_i}(w) =\alpha_i(u).$ We thus
have associated to $w\in {\tilCrit}(\gt^*,z)$ an element $u\in \ga$.
If none of the coordinates $z_i$ of $z$ vanishes, then
$\alpha_i(u)\neq0$ for $\dt in$, and, consequently, $u\in U(\A)$.

Now consider the point $u$ associated to $w$ this way, and let
$\lambda=\sum_{i=1}^n n_i \omega_i$ be an element in $\ga_\Z$, i.e.
$\sum_{i=1}^n n_i \mu_i\in \R \gamm$.  Then
\[ \prod_{i=1}^n \alpha_i(u)^{n_i}=\prod_{i=1}^n
z_i^{n_i}e_ {\sum_{i=1}^n n_i \beta_i}(w)=z^{\lambda},\] since
$w\in\Torus(\gt^*)$.  Thus we obtain that the equation
$p_{\lambda}(u)=z^{\lambda}$ holds for all $\lambda\in \ga_{\Z}$,
i.e. $u\in O(z,\A)$.

Conversely, if $u\in O(z,\A)$, then, as shown in
(\cite{tor1},\,Lemma 4.3), we can find $w\in \Torus(\gt^*)$ such
that $z_i e_{ \beta_i}(w) =\alpha_i(u)$; this element $w$ is  in
$\tilCrit(\gt^*,z)$.
\end{proof}

Now we turn to \eqref{localform}; our next step is to express the
values $f_z(w)$ and $\Hesszt(w)$ for $w\in\tilCrit(\gt^*,z)$ as values
of certain functions on $\ga_{\C}$ evaluated at the corresponding
point $u$.  Introduce the polynomial $G$ on $\ga_\C$:
\begin{equation}\label{G}
G(u)=\sum_{\nu}\vol_{{\B}}(\nu)^2\prod_{j\notin
\nu}\alpha_j(u),
\end{equation} where $\nu$ runs over the  subsets of cardinality $d$
of the set $\{1,2,\ldots, n\}$. In this equation the volume
$\vol_{{\B}}(\nu)$ is the volume of the parallelepiped $\sum_{j\in
  \nu} [0,1]{\beta_j}$ in $\gt$. In particular, in the sum in \eqref{G} only
those subsets $\nu$ will contribute for which the vectors
$\{\beta_j,j\in\nu\}$ are linearly independent. Incidentally, these
are exactly the subsets for which the vectors
$\{\alpha_j,\,j\notin\nu\}$ are linearly independent (cf.
\cite[Proposition 1.2]{tor1}).

We will denote by $w(u)$ the point of $\tilCrit(\gt^*,z)$
associated to a point $u\in O(z,\A)$, i.e.
 when $w$ and
$u$ are related by \eqref{w}.

\begin{lemma}\label{formulaG}
  If $w=w(u)$, then
$$f_z(w(u))=\sum_{i=1}^n \alpha_i(u),$$
$$\Hesszt(w(u))=G(u).$$
\end{lemma}
\begin{proof} The first formula follows immediately from the definition of
$f_z=\sum_{i=1}^n z_i e_{\mu_i }$  and equation (\ref{w}).
The second formula is  proved in \cite[Proposition 4.5]{tor1}.
\end{proof}

We introduce the notation $\kappa_\A=\sum_{i=1}^n \alpha_i\in \ga^*$.

\begin{prop}\label{local}
  Let $P$ be a homogeneous polynomial of degree $d-\ss+1$ in $n$
  variables; then $e_{\gamm}P(z_1 e_{\mu_1},\ldots,
  z_n e_{\mu_n})$ is an element of $I^{d+1}(C)$.  If $\gamma\in
  V^*_\Z$ and $z\in \C^n$
  satisfy conditions \eqref{cond1}, \eqref{cond2} and \eqref{cond3}, then
  \begin{equation}
    \label{torlocal}
\tres_{\bff}\, e_{\gamm}P(z_1 e_{\mu_1},\ldots,
  z_n e_{\mu_n})={\ss} \sum_{u\in O(z,\A)}
\frac{P(\alpha_1(u),\ldots,\alpha_n(u))}{\kappa_\A(u) G(u)}.
  \end{equation}
\end{prop}
\begin{proof}
Using Proposition \ref{tlocal} and \eqref{obvious}, we may write
\begin{multline*}
\tres_{\bff}\, e_{\gamm}P(z_1 e_{\mu_1},\ldots,
  z_n e_{\mu_n})=\\
\ss\sum_{w\in \tilCrit(\gt^*,z)} \frac{e_{\gamm}(w)
  P(z_1 e_{\beta_1}(w),\ldots, z_n e_{\beta_n}(w))}{ f_z(w)
  \Hesszt(w)}.
\end{multline*}
Note that the character $e_\gamm$ is identically $1$ on
$\Torus(\gt^*)$, thus Lemma \ref{formulaG} implies the Proposition.
\end{proof}

\section{Summary of our earlier results}
\label{sec:summary}

\subsection{Flags and iterated residues}

Our goal in this section is to develop a Laurent-type expansion for
the functions of the complex variable $z=(z_1,\dots,z_n)$, which
appear on the right hand side of \eqref{torlocal}. These are
essentially results of \cite{tor1}, but the exposition is improved
here. In particular, we present an invariant formalism, which avoids
the choice of a special ($\gc$-positive) basis of $\ga_\Z$.  We also
introduce a more precise notion of regularity of our target vector
$\xi$.

We maintain some of the notation of the previous section. Thus we have an
exact sequence
$$0\lra \ga\overset{\iota}{\lra}\gg\overset{\pi}\lra  \gt\lra 0$$
of
real vector spaces of dimensions $r$, $n$ and $d$; thus $n-d=r$. The
space $\gg$ has a fixed basis $(\omega_1,\dots,\omega_n)$, and we have
three lattices
\[\ggz=\oplus_{i=1}^n\Z\omega_i,\; \gtz=\pi(\ggz)\;,
\gaz=\ker(\pi|_{\ggz}).
\]
This construction gives rise to the Gale dual sequences of vectors
\[ \A=[\alpha_1,\dots,\alpha_n]\subset\ga_\Z^*,\quad
\B=[\beta_1,\dots,\beta_n]\subset\gt_\Z.\]

Denote by $\rata$ the graded linear space of rational functions on
$\ga_\C$ whose denominators are products of powers of elements of
$\A$. The functions in $\rata$ are regular on $U(\A)$, the complement
of the hyperplane arrangement induced by $\A$. Recall that for each
$\lambda\in \ga_\Z$, we defined a rational function
$p_\lambda\in\rata$ in \eqref{plam} of homogeneous degree
$\langle\kappa_\A,\lambda\rangle$, where $\kappa_\A=\sumn\alpha_i$.

We fix an orientation of $\ga$, and choose an oriented $\Z$-basis
$\bl=[\lambda_1,\ldots,\lambda_r]$ of $\gaz$. We denote by $\dvol$
the holomorphic $r$-form $d\gamma_1\wedge d\gamma_2\wedge \cdots
\wedge d\gamma_r$ on $\ga_\C$, where
$[\gamma_1,\dots,\gamma_r]\subset\ga^*$ is the basis dual to
$\bl$.  This form depends only on the orientation and the integral
structure of $\ga$, and not on the choice of the basis $\bl$.

Introduce the notation $p_j=p_{\lambda_j}$ and $q_j=z^{\lambda_j}$,
and define the maps
\begin{equation} \label{pbl}
p_{\bl}=(p_1,p_2,\ldots, p_r):U(\A)\ar(\C^*)^r
\end{equation}
and
$$q_{\bl}=(q_1,q_2,\ldots,q_r):(\C^*)^n\ar(\C^*)^r.$$

The map $p_\bl$ plays a central role in our investigations. We will
need a formula (\cite[Proposition 4.5]{tor1}) for the Jacobian of this
map:
\begin{lemma}\label{thejac}
  We have
  the equality of meromorphic $r$ forms on $U(\A)$:
\begin{equation}\label{Jac}
\frac{dp_1}{p_1}\wedge\frac{dp_2}{p_2}\wedge\cdots\wedge\frac{dp_r}{p_r}=
\frac{G(u)\,\dvol}{ \prod_{i=1}^n\alpha_i(u)},
\end{equation}
where the function $G(u)$ is defined in \eqref{G}.
 \end{lemma}

 There is, in fact, a more invariant way to describe the map
 $p_{\bl}$.  Let $p:\comp\ar\Torus(\ga^*)$ be the equidimensional map given by the
 formula
\[ p(u) = \frac{1}{2\pi\sqrt{-1}}\sumn \log(\alpha_i(u))\,\alpha_i.\]
The map formally depends on the branch of the logarithm one chooses,
however, this choice is immaterial in the quotient
$\Torus(\ga^*)=\ga^*_{\C}/\ga^*_{\Z}$.

Similarly let $q:(\C^*)^n\ar\Torus(\ga^*)$ be the map given by the
 formula
\[ q(z) = \frac{1}{2\pi\sqrt{-1}}\sumn \log(z_i)\,\alpha_i.\]

If $\lambda\in \ga_\Z$, then $e_\lambda(p(u))=p_\lambda(u)$ and
$e_\lambda(q(z))=z^{\lambda}$. Thus the set $O(z,\A)$ (see
Definition \ref{defo}) may be given as the set of  solutions of
the equation on $U(\A)$:
\[ O(z,\A)= \{u\in \comp;\, p(u)=q(z)\};\]
one can also write $O(z,\A)=p^{-1}(q(z))$.

Using the $\Z$-basis $\bl=[\lambda_1,\ldots,\lambda_r]$ of $\gaz$, we
identify $\Torus(\ga^*)$ with $(\C^*)^r$ via the map $ (e_{\lambda_1},
\ldots, e_{\lambda_r})$. After this identification, the map $p_{\bl}$
is the map $p$ and $q_{\bl}$ is the map $q$. We collected the maps we
will use in the following diagram:
\begin{equation}
  \label{diagram2}
\begin{diagram}
  \node[2]{\comp}\arrow{sw,t}{p_{\bl}}\arrow[2]{s,r}{\iota}
\arrow{se,b}{\xia}\arrow{ese,t}{p}\\
  \node{(\C^*)^r} \node[2]{\ga^*}
  \node{\Torus(\ga^*)}\arrow{w,b}{\II}\\
  \node[2]{(\C^*)^n}\arrow{nw,b}{q_{\bl}}\arrow{ne,t}{\xig}
  \arrow{ene,b}{q}
\end{diagram}
  \end{equation}

\begin{itemize}
\item
$\iota$  The map from $\comp$  to $(\C^*)^n$ is given by $u\mapsto
(\alpha_1(u),\ldots,\alpha_n(u))$; it is the restriction of the map
$\iota$ in \eqref{exact1};
\item $\II$ is the map $w\mapsto -\mathrm{Re}(2\pi\sqrt{-1}w)$ from
$\ga^*_{\C}/\ga^*_{\Z}$ to $\ga^*$;
\item $\xig:(z_1,\dots,z_n)\mapsto -\sumn\log|z_i|\,\alpha_i$.
\item
  $\xia(u)=\xig(\iota(u))=\II(p(u))=-\sumn\log(|\alpha_i(u)|)\,\alpha_i(u)$.
\end{itemize}
The following statement is straightforward:
\begin{lemma}\label{easyproper}
The maps $\iota$ and $\II$ are proper.
\end{lemma}

To formulate our expansion formula, we need some preparations.  We
introduce the following notions from \cite{tor1}.  Let $\flag(\A)$ be
the finite set of flags
\[F = [F_0=\{0\}\subset F_1\subset F_2\subset  \dots \subset F_{r-1}
\subset F_r=\ga^*],\text{ where }\dim  F_j=j,\] such that the
finite sequence $\A\cap F_j$ spans $F_j$ for each $j=1,\dots,r$.
For each $F\in \flag(\A)$, we choose an ordered basis
$(\gamma^F_1,\gamma^F_2,\ldots, \gamma^F_r)$ of $\ga^*$, which is
unimodular with respect to the volume form $\dvol$ on $\ga$, and which
is such that $F_j=\oplus_{k=1}^j \R \gamma^F_k$. If $F$ is fixed, then
we will use the simplified notation $u_j=\langle \gamma_j^F,u\rangle$,
$\dt jr$. We call the functionals $u_j$, $\dt jr$, coordinates {\em
  adapted to} $F$.

Now let $F\in \flag(\A)$, and let $N$ be a positive real number.
Define the open subset $U(F,N)$ of $\ga_\C$ by
\[
U(F, N)=\{u\in \ga_\C\,;\; 0< N|u_j|< |u_{j+1}|,\, j=1,2,\ldots,
r-1\}.
\]

We make the following simple observations:
\begin{lemma}\label{ufn}
  \begin{enumerate}
  \item If $N$ is sufficiently large, then the sets
    $\{U(F,N);\;F\in\flaga\}$ are disjoint, moreover, for each
    $F\in\flaga$, $U(F,N)\subset U(\A)$, and the linear form
    $\kappa_\A=\sum_{i=1}^n\alpha_i$ does not vanish on $U(F,N)$.
\item  The set $U(F,N)$ is
invariant under the scaling $u\to e^t u$.
\item If $N_2>N_1$, then
$U(F,N_2)\subset U(F,N_1)$.
\item For $N$ sufficiently large, $U(F,N)$ is isomorphic to the product
  of $\C^*$ with $r-1$ small punctured disks. Choose a sequence
  of real numbers $\ve:0<\epsilon_1\ll\epsilon_2\ll\dots
  \ll\epsilon_r$, where $\epsilon\ll\delta$ means
  $N\epsilon<\delta$. Then the embedded torus
\begin{equation}
  \label{zefepsilon}
T_F(\ve) = \{u\in\ga_{\C};\; |u_j|=\epsilon_j,\,
j=1,\dots,r\}\subset U(F,N)
\end{equation}
oriented by the form $d\arg u_1\wedge\dots\wedge d\arg u_r$,
represents a generator of $H_r(U(F,N),\Z)$.
\end{enumerate}
\end{lemma}
\begin{definition}\label{hf}
  Let $F\in\flaga$, and choose $N$ large enough in order to have
  $U(F,N)\subset\comp$. We denote by $h(F)$ the homology class in
  $H_r(\comp,\Z)$ of the oriented compact torus $T_F(\ve)$.
\end{definition}
If $\omega$ is a holomorphic $r$-form on $U(F,N)$, then it is closed.
Thus we can write $\int_{h(F)}\omega$ for the integral
$\int_{T_F(\ve)}\omega$, since it only depends on the homology
class $h(F)$ of the cycle $T_F(\ve)$.

If $\omega=\phi(u)\dvol$, where $\phi\in \rata$, then
$\frac{1}{(2\pi\sqrt{-1})^r}\int_{h(F)}\phi(u)\dvol$ coincides with
the {\em iterated residue} of $\phi$ with the respect to $F$
(\cite{asz_ir}). This is defined as follows: write $\phi$ as a rational
function $\phi^F$ of the coordinates $u_j$ and define the iterated
residue as
\begin{equation}\label{res}
\res_F(\phi)=\res_{u_r=0}\,du_r\;\res_{u_{r-1}=0}\cdots
\res_{u_1=0}\,du_1\;\phi^F(u_1,u_2,\ldots,u_r),
\end{equation}
where each residue is taken assuming that the
variables with higher indices have a fixed, nonzero value.

Let $F\in \flag(\A)$ be a flag.  Introduce the vectors
\[
\kappa^F_j=\sum\{\alpha_i;\;i=1,\dots,n,\, \alpha_i\in F_j\},\;
j=1,\ldots, r.
\]
Note that $\kappa^F_r$ is independent of $F$ and equals
$\kappa_\A=\sum_{i=1}^n\alpha_i$.

\begin{definition}\label{flags}
  \begin{itemize}
  \item We say that a flag $F\in \flag(\A)$ is {\em proper} if the vectors
$\kappa_j^F$, $j=1,\ldots, r$, are linearly independent.

\item Assuming $F$ is proper, define the sign
$\nu(F) = \pm 1$ depending on whether the sequence
$[\kappa_1^F,\ldots, \kappa_r^F]$ is positively oriented or not.
\item For a flag $F\in \flaga$, introduce the non-acute cone $\s(F,\A)$
generated by the non-negative linear combinations of the elements
$\{\kappa^F_j, j=1,\ldots,r-1\}$ and the line $\R \kappa_\A$:
$$\s(F,\A) =\sum_{j=1}^{r-1}\R^{\geq 0}\kappa_j^F+\R \kappa_\A$$
\item For $\xi\in \ga^*$ denote by $\FF$ the set of flags $F\in
\flaga$ such that $\xi\in\s(F,\A)$.
  \end{itemize}
\end{definition}
\begin{remark} Note that for any $t\in \R$, we have
$\flag(\A,\xi)=\flag(\A,\xi+t\kappa_\A)$.
\end{remark}

In \cite{tor1}, we introduced the notion of $\tau$-regularity for
a vector $\xi\in\ga^*$. In the present paper we replace it with a
less restrictive notion of being removed from the boundary inside
the cones $\s(F,A)$, $F\in\FF$.

\begin{definition}\label{regular}
  We define a vector $\xi\in\ga^*$ to be {\em $\flaga$-regular} if for
  any flag $F\in\flag(\A,\xi)$ the first $r-1$ coefficients,
  $m_1,\dots,m_{r-1}$, in any linear expression
 \[ \xi = m_1\kappa^F_1+\cdots +m_{r-1}\kappa^F_{r-1}+m_r\kappa^F_r \]
 do not vanish. Given a positive number $\tau$, we will call $\xi$
 {\em $\flagat$-regular} if in any such linear expression we
 have $\min_{1\leq j\leq r-1} m_j>\tau$. We denote by $\tra$ the
 subset  of $\flagat$-regular elements of $\ga^*$.
 \end{definition}
 \begin{remark}
   \begin{enumerate}
   \item Note that in the definition, $F\in\flag(\A,\xi)$ implies that
     $m_j\geq0$, $\dt jr$
\item It is easy to see that if $\xi$ is a $\flaga$-regular
 element of $\ga^*$, then every flag $F\in \FF$ is proper.
\item If $\xi$ is a $\flagat$-regular, then for any $t\in \R$,
  the element $\xi+t\kappa_\A$ is also $\flagat$-regular.
\item  The set $\tra$ is a disjoint
  union of open cones in $\ga^*$. It occupies ``most'' of the
  vector space $\ga^*$, in the sense that the intersection of its complement
with a generic affine line $L\subset\ga^*$ is bounded.
   \end{enumerate}
\end{remark}

 Consider the maps on the diagram
\eqref{diagram2}.  For $\xi\in \ga^*$,  define the set
\begin{equation}
  \label{defzxi}
Z(\xi)=\left\{u\in\comp;\;
\sumn\log|\alpha_i(u)|\alpha_i=-\xi\right\}
\end{equation}
Observe that
\[ Z(\xi) = \xia^{-1}(\xi)=(\II\circ p)^{-1}(\xi)=
(\xig\circ \iota)^{-1}(\xi).\]

Choosing an oriented $\Z$-basis $\bl=[\lambda_1,\dots,\lambda_r]$ of
$\ga_\Z$, we can present the set $Z(\xi)$ as the set of common
solutions of the $r$ analytic equations on $U(\A)$:
$$Z(\xi)=\{|p_j(u)|=e^{-\langle \xi,\lambda_j\rangle} ;\;j=1,\dots,
r\};$$
thus $Z(\xi)$ is the inverse image of the product of $r$
circles under the map $p: U(\A)\to \Torus(\ga^*)=(\C^*)^r$.  In
particular, $Z(\xi)$ is a real analytic subvariety of $\comp$.  The
form $d\arg p_1\wedge\dots \wedge d\arg p_r$ defines an orientation of
$Z(\xi)$; this orientation depends only on the orientation of $\ga$
and not on the particular basis oriented $\bl$ we picked.

\begin{example}
\label{ex1}
  Consider the case $d=r=2, n=4$, $\ga=\R^2$, $\gaz=\Z^2$, with
  coordinates $x,y$. Let $\A=(y,x+y,y,x)$. Then $\kappa=2x+3y$, and
  $\kappa$ is in the interior of the chamber $\{mx+ny;\;0<m<n\}$,
  which corresponds to the Hirzebruch surface, the blow-up of
  $\mathbb{P}^2$ at one point. Pick the vector
  $\xi=(\log\epsilon)(x+2y)$ in this chamber, where $\epsilon$ is a
  small positive constant. Then the equations defining $Z(\xi)$ read as
\begin{equation}
  \label{zex}
|x(x+y)|=\epsilon,\quad |y^2(x+y)|=\epsilon^2.  
\end{equation}
\qed
\end{example}

With these preparations, we are ready to formulate the main results of
\cite{tor1}.  They were proved assuming the notion of
$\tau$-regularity for a vector $\xi\in\ga^*$ in the strong sense of
Definition 2.2 of \cite{tor1}. It is easy to see, however, that the
results, as well as the proofs of \cite{tor1} remain true under the
weaker assumption of $\flagat$-regularity.

\begin{definition}
  We say that a continuous map $\tilde p:U\ar V$ is {\em proper to} an
  open subset $V'\subset V$ if for every compact $K\subset V'$ the
  inverse image $\tilde{p}^{-1}(K)$ is compact.
\end{definition}

\begin{theorem}\label{main1}
  \begin{enumerate}
  \item For any sufficiently large $N>0$ and any proper flag
    $F\in\flaga$, the holomorphic map $p$, restricted to $U(F,N)$ is
    non-singular (\cite[\S5.3]{tor1}).
\item For $\tau$ sufficiently large, the map $\xia:U(\A)\ar\ga^*$
  is proper to the set $\tra$ of $\flagat$-regular vectors in
  $\ga^*$.
\item Given $N>0$, there is $\tau>0$ such that if $\xi\in\tra$ then
\begin{equation}
  \label{Zin}
Z(\xi) \subset\cup_{F\in\FF}U(F,N).
\end{equation}
\end{enumerate}
\end{theorem}
It will be convenient to make the following
\begin{definition}
 Fix $N$ sufficiently large to satisfy the conditions of Lemma
 \ref{ufn} (1) and Theorem \ref{main1} (1), and let $\tau$ be such that it
 satisfy the conditions of Theorem \ref{main1} (2) and also
 statement (3) of the same Theorem with respect to $N$. Then we will call
 the pair of positive constants $(N,\tau)$ {\em sufficient}.
\end{definition}
Let us point out some corollaries of the Theorem.
\begin{cor}\label{impcor}
  \begin{enumerate}
  \item The map $p$ is proper to the set $\II^{-1}(\tra)$. Similarly,
  the map $p_\bl$ is proper to the set of vectors of the form
\[\left\{(q_1,\dots,q_r);\;
  |q_j|=e^{-\lr\xi{\lambda_j}},\,\dt jr\,
\text{ for some }\xi\in\tra\right\}.
\]
\item For sufficient $(\tau,N)$ and any $z\in (\C^*)^n$ such that
  $\xig(z)\in\tra$, the set $O(z,\A)$ is finite, and has the
  following decomposition into a disjoint union:
$$O(z,\A)=\bigcup
\left\{O(z,\A)\cap U(F,N);\;F\in \flag(\A,\psi(z))\right\}.$$
\item For sufficient $(N,\tau)$ and any $\xi\in \tra$,
the cycle $Z(\xi)$ is smooth and
  compact in $\comp$.
  \end{enumerate}
\end{cor}
The first statement follows from Lemma \ref{easyproper} and Theorem
\ref{main1}.  To prove the second statement, recall that
\begin{itemize}
\item $O(z,\A)=p^{-1}(q(z))$,
\item for $\xi\in\tra$ every flag in $\flag(\A,\xi)$ is proper,
  and
\item $\displaystyle \II^{-1}(\tra) = q\left(\psi^{-1}(\tra)\right).$
\end{itemize}
 Thus we represented $O(z,\A)$ as the inverse
image of a point under a proper non-singular equidimensional map. This
implies that $O(z,\A)$ is finite. The second part of statement (2)
follows from Theorem \ref{main1} (3).

Finally, Statement (3) of Corollary \ref{impcor} follows from a
similar argument, since $Z(\xi)$ is the inverse image of a smooth
torus under $p$.
\begin{remark}\label{importantrem}
   We would like to emphasize that, in general, the map $p$ is {\em
    not} proper on the whole of $\comp$, and this is why the homology
  class of the cycle $Z(\xi)$ may vary with $\xi$.
\end{remark}

Once we know that the cycle $Z(\xi)$ is compact and
smooth in $\comp$, it is natural to try to compute its homology class
in $H_r(\comp,\Z)$. The method of this computation relies on a
degeneration technique, which is known as {\em tropical geometry}
\cite{Viro,Sturm}. Let us explain this on our Example \ref{ex1}.
We consider \eqref{zex}, let the parameter $\epsilon$ be sufficiently
small, and introduce the tropical ansatz 
\[  |x|=\epsilon^a,\;|y|=\epsilon^b,\;|x+y|=\epsilon^c.
\]
Then the equations \eqref{zex} imply the equalities
\begin{equation}
\label{system}
  \begin{cases}
    a+c&=1,\\ 2b+c&=2.
  \end{cases}
\end{equation}
These, naturally, do not determine $a,b$ and $c$, but we observe that
for small enough $\epsilon$, if $|x+y|$ is very small compared to $|x|$,
then $|x|$ and $|y|$ should be rather close.  These gives us the following
three possibilities:
\begin{enumerate}
\item $c>a,b$, which implies $a\sim b$
\item $b>a,c$, which implies $a\sim c$
\item $a>c,b$, which implies $c\sim b$.
\end{enumerate}
Here by $>$, we mean {\em significantly greater}, and by $\sim$ we
mean {\em very close}.  Now we can go back to our system
\eqref{system}, and solve them under the three possible conditions
$a=b$, $a=c$, or $b=c$, and obtain the three solutions
\[ (1,1,0),\quad(1/2,3/4,1/2),\quad(1/3,2/3,2/3).
\]
However, only the second of these equations satisfies the
corresponding inequalities. Thus we, informally, conclude, that for
$\xi=x+2y$ the cycle $Z_\epsilon(\xi)$ consists of a torus which is
very close to the torus
$\{|y|=\epsilon^{3/4},\,|x|=\epsilon^{1/2}\}\subset\C^2$. Integration
over such a torus is equivalent to a single iterated residue:
\[ \JK_{\c(\kappa)}(f) = \res_x\res_y f\, dx\,dy.
\]

The result in the general case is the following

\begin{theorem}\label{homclass}
  For a sufficient pair $(N,\tau)$ and $\xi\in\tra$, the homology class
  $[Z(\xi)\cap U(F,N)]\in H_r(\comp,\Z)$ is equal to $\nu(F)h(F)$,
(cf. Definitions \ref{flags},\ref{hf}).  Hence, using Theorem
  \ref{main1} (3), we can conclude that
  \begin{equation}
    \label{zhom}
[Z(\xi)] = \sum_{F\in\FF}\nu(F)h(F)\in H_r(\comp,\Z),
  \end{equation}
\end{theorem}

Recall that Lemma \ref{OZA} establishes a bijection $w\rightarrow
u$ between $\tilCrit(z)$ and $O(z,\A)$, and that under this
bijection the value of the function $f_z\Hesszt$ at $w$ coincides
with the value of the function $\kappa_\A G$ at $u$ (Lemma
\ref{formulaG}). The function $\kappa_\A(u) G(u)$ does not vanish
on $U(F,N)$ provided $N$ is sufficiently large and $F$ is a proper
flag. The fact that $\kappa_\A$ does not vanish on $U(F,N)$ for
large $N$  is easy and stated in Lemma \ref{ufn}. The fact that
$G(u)$ does not vanish on $U(F,N)$ provided $N$ is sufficiently
large and $F$ is  proper is the content of Proposition 5.9 of
\cite{tor1}. Thus we obtain the following.

\begin{cor}\label{verifycond3}
  Let $z\in (\C^*)^n$ such that $\xig(z)$ is
  $\flagat$-regular.  Then $\tilCrit(z)$ is finite and the
  function $f_z\Hesszt$ does not vanish at any of the points of
  $\tilCrit(z)\subset\Toric(C)$.
\end{cor}

\begin{prop}\label{sum} Let $(\tau,N)$ be a sufficient  pair of
  constants, and let $z=(z_1,z_2,\ldots, z_n)\in(\C^*)^n$ be such that
  $\xig(z)$ is $\flagat$-regular. Then for any $F\in
  \flag(\A,\xig(z))$ and any holomorphic function $\phi$ on
  $U(F,N)$, we have
  \begin{equation}
    \label{intheprop}
\sum_{u\in O(z,\A)\cap U(F,N)}\phi(u)=\frac{\nu(F)}{(2\pi\sqrt{-1})^r}
\sum_{\lambda\in \ga_\Z}
\int_{h(F)}\frac{\phi}{p_\lambda} \frac{dp_1}{p_1}\wedge
\cdots\wedge \frac{dp_r}{p_r}z^{\lambda},
  \end{equation}
  and the sum is absolutely convergent on the domain $\psi^{-1}(\tra)$.
\end{prop}

\begin{proof} We need  two facts from complex function theory. First
  recall the Laurent expansion of a function of one complex variable.
  Let $\phi$ be an holomorphic function on an annulus
  $\epsilon_1<|z|<\epsilon_2$.  Then, for
  $\epsilon_1<\epsilon,|q|<\epsilon_2$, we have
  \begin{equation}
    \label{laurant}
\phi(q)= \frac{1}{2\pi\sqrt{-1}}\sum_{n\in \Z} \int_{|y|=\epsilon}
  \frac{\phi(y)}{y^n}\frac{dy}{y}q^n,
  \end{equation}
  and the series is uniformly and
  absolutely convergent on any compact subset of $
  \epsilon_1<|q|<\epsilon_2$.  A similar statement holds for functions
  of $r$ complex variables, which are defined on a product on annuli.

The other fact is contained in the following
  \begin{lemma}\label{pushfor}
  Let $\tilde p:U\ar V$ be a proper holomorphic map between open
  subsets of a complex vector space, and let $\phi:U\ar\C$ be a
  holomorphic function. Then the push-forward function $\tilde
  p_*\phi$ given by
\[ \tilde p_*\phi(v)= \sum_{p(u)=v} \phi(u) \]
is holomorphic on $V$.
      \end{lemma}
The proof will be omitted (cf \cite[Chapter 5]{GH}).
We will, in fact, apply this to our map $p$ in a domain
where it is proper and {\em non-singular}; in this case the statement
is trivial.

Now we are ready to prove the Proposition. Let $\xi\in\ga^*$, and
orient the torus
\[
\TT=\{q\in(\C^*)^r;\;|q_j|=
e^{-\langle\xi,\lambda_j\rangle},\,j=1,\dots,r\}\subset(\C^*)^r
\]
in the standard fashion. According to Corollary \ref{impcor} (1), if
$\xi\in\tra$, then the map $p_{\bl}=(p_1,\dots,p_r)$, even when
restricted to $U(F,N)$, is proper to a neighborhood of $\TT$.

Now set $\xi=\psi(z)$, and consider the function $\phi$ given in
the proposition. From the discussion above, we can conclude that
the function $p_{\bl*}\phi$ is holomorphic in a neighborhood of
the torus $\TT$, and thus we can write down the standard Laurent
expansion for it:

\[\sum_{n\in\Z^r}\frac{1}{(2\pi\sqrt{-1})^r}
\int_{\TT}\frac{p_{\bl*}\phi}{y_1^{n_1}\cdots y_r^{n_r}}
\frac{dy_1}{y_1} \frac{dy_2}{y_2}\ldots
\frac{dy_r}{y_r}q_1^{n_1}\ldots q_r^{n_r}.
\]

Now we pull back this equality by the map $p_{\bl}$, restricted to
$U(F,N)$. For the left hand side we have
\[ p_{\bl*}\phi(q_\bl(z)) = \sum\phi(u),\;u\in O(z,\A)\cap U(F,N).\]
To compute the pull-back of the right hand side, observe that
$p_\bl^{-1}(\TT)=Z(\xi)$, thus using Theorem \ref{zhom} we can
conclude that this pull-back equals
$$\sum_{n\in\Z^r}\frac{\nu(F)}{(2\pi\sqrt{-1})^r}
\int_{h(F)}\frac{\phi}{p_1^{n_1}\cdots p_r^{n_r}} \frac{dp_1}{p_1}
\frac{dp_2}{p_2}\ldots \frac{dp_r}{p_r}q_1^{n_1}\ldots
q_r^{n_r}.$$

Thus we recovered the two sides of \eqref{intheprop}, and this
completes the proof of the proposition.
\end{proof}

Combining Proposition \ref{sum} with Lemma \ref{thejac}, Corollary
\ref{impcor} and Theorem \ref{homclass}, we obtain
\begin{cor}\label{FORMULA}
Assume that $(\tau,N)$  is a sufficient pair of constants.  Let
$z\in (\C^*)^n$ be such that $\xig(z)$ is
$\flagat$-regular.
 If $K$ is a holomorphic function defined
 on $\cup_{F\in \flag(\A,\xig(z))} U(F,N)$, then we have:
$$\sum_{u\in O(z,\A)} \frac{K(u)}{G(u)}=
\frac{1}{(2\pi\sqrt{-1})^r}\sum_{F\in \flag(\A,\xi)}\nu(F)\sum_{\lambda\in
\ga_\Z}\int_{h(F)}
\frac{K\,\dvol}{p_\lambda\prod_{i=1}^n\alpha_i} z^{\lambda},$$ where
$G$ is defined in \eqref{G}. The series is
absolutely convergent on $\xig^{-1}(\tra)$.
\end{cor}

\subsection{Jeffrey-Kirwan residue and iterated residues}

Consider our  exact sequence
\begin{equation}
\label{es}
 0\ra\ga\rightarrow\gd\overset\pi\rightarrow\gt\ra0
\end{equation}
and the dual sequence
\begin{equation} \label{des}
 0\ra\dgt\rightarrow\dgd\overset\themu\rightarrow\dga\ra0
\end{equation}
which is is also exact.

The sequence $\CB$ is the image under $\pi$ of the canonical basis
$\{\omega_i\}_{i=1}^n$ of $\gd$, while the sequence $\CA$ is the image under
$\themu$ of the dual canonical basis $\{\omega^i\}_{i=1}^n$ of $\dgd$. In the
toric mirror residue conjecture of Batyrev-Materov \cite{BM}, two
toric varieties make their appearance. The first one is the toric
variety associated to the polytope $\Pi^{\B}\subset \gt$ obtained
as the convex hull of $\CB$. The other one $ V_\A(\gc)$ is
determined by a chamber $\gc$ of the cone generated by $\A$
containing $\kappa_\A$ in their closure.

Recall the definition of a chamber. Consider the cone
$\Cone(\A)\subset\ga^*$ generated by the elements of $\A$. The
subset $\Cone_{\sing}(\A)$ is the set of elements in $\Cone(\A)$
which can be written as a positive linear combination of $m$
elements of $\A$, with $m<r$. Chambers are the open polyhedral
cones in $\ga^*$ which are the connected components of
$\Cone(\A)\setminus \Cone_{\sing}(\A)$. An element $\xi$  in
$\Cone(\A)\setminus \Cone_{\sing}(\A)$ will be called
$\A$-regular; it belongs then to a unique chamber $\gc$.

Now we recall the relation between Minkowski sums of polytopes and
chambers.

\begin{definition}\label{partition} Let $\theta\in \ga^*$ and define
the
 {\em
  partition polytope}
 \[
\Pi_\theta=\themu^{-1}(\theta)\bigcap\sum_{i=1}^n
\R^{\geq0}\omega^i,
\]
which lies in an affine subspace of $\gd^*$ parallel to $\gt^*$.
\end{definition}
The following proposition is well-known; see \cite[Lemma 3.4]{poly}
for a proof.

\begin{prop}\label{minkowski}
Let $\theta_k$, $k=1,\ldots, l$ be elements in $\ga^*$. Let
$\theta=\sum_{k=1}^l\theta_k$. Then  the partition polytope
$\Pi_{\theta}$ is the Minkowski sum of the partition polytopes
$\Pi_{\theta_k}$ if and only the elements $\theta_k$ are in the
closure of the same chamber $\gc$.
\end{prop}

A chamber $\gc$ determines a complete regular simplicial fan in $\gt$;
the $1$-dimensional faces of these fans are among the rays
$\R^+\beta_i, \beta_i\in \CB$.  A cone $\sum_{i\in
  \nu}\R^+\beta_i\subset \gt$ belongs to the fan determined by $\gc$
if and only if the cone $\sum_{k\notin \nu}\R^+\alpha_k$ contains the
chamber $\gc$. Denote the toric variety associated to this fan by
$V_{\A}(\gc)$. It is an orbifold and it can also be realized by a
quotient construction (cf \cite{tor1}, Section 1.2). In particular,
each polynomial function $P$ on $\ga$ leads to a cohomology class
$\chi(P)$ on $V_\A(\gc)$ via the Chern-Weil map.

Now we recall a formula from \cite{tor1} for the intersection
numbers on $V_\A(\gc)$ for any chamber $\gc$.  At the end of this
section, the special chambers containing $\kappa_\A$ in their
closure will appear.

Associated to each chamber $\gc$, there is a linear form $\phi \mapsto
\JK_{\gc}(\phi)$ on $\rata$ called the {\em Jeffrey-Kirwan
  residue}. This linear form vanishes on homogeneous elements of $\rata$
unless its homogeneous degree is equal to $-r$.  We refer to
\cite{B-V} for the definition. Here we only recall an important
observation to be used later:

 Let $\gcp\subset \ga$ be the polar cone of the chamber
 $\gc$:
 \[  \gcp=\{\lambda\in\ga;\; \langle  \xi,\lambda\rangle  \geq 0,
\text{ for all } \xi\in\gc\}. \]

\begin{lemma}{\cite[Proposition 3.1]{tor1}}\label{null}
  Let $\lambda\in \ga_Z$ such that $\lambda$ does not belong to $\gcp$,
  and let $P$ be a polynomial on $\ga$. Then
$$\JK_{\gc}\left(\frac{P}{p_\lambda\prod_{i=1}^n\alpha_i}\right)=0$$
\end{lemma}

As explained in \cite{tor1}, an algebraic  formula for integrating
a cohomology class over  $V_\A(\gc)$ can be given in terms of the
Jeffrey-Kirwan residue.  For a polynomial $P$ on $\ga$:
\begin{equation}
  \label{interjk}
\int_{V_\A(\gc)}\chi(P)=JK_{\gc}\left(\frac{P}
  {\prod_{i=1}^n\alpha_i}\right).
\end{equation}

We recall now our construction of a concrete, compact, real-analytic
cycle in $\comp$, such that the linear form $\JK_{\gc}(\phi)$ is simply
obtained by integration of $\phi(u)\,\dvol$ on this cycle.

\begin{definition}\label{defsplus}
For a flag $F\in \flaga$, introduce the acute cone $\s^+(F,\A)$
generated by the non-negative linear combinations of the elements
$\{\kappa^F_j, j=1,\ldots,r\}$
 \[ \s^+(F,\A) = \sum_{j=1}^r \R^{\geq 0}\kappa^F_j.\]
For $\xi\in \ga^*$ denote by  $\FFF$ those flags in
$\flaga$ for which  $\xi\in\s^+(F,\A)$.

We say  that $\xi\in\ga^*$ is {\em $\flaga^+$-regular} if $\xi$ is
not on the boundary of any of the cones $\s^+(F,\A)$, $F\in\FFF$.
This means that for any flag $F\in\FFF$, all the coefficients,
$m_1,m_2,\dots,m_r$, in any linear expression
  \[ \xi = m_1\kappa^F_1+m_2\kappa^F_2+\cdots
  +m_{r-1}\kappa^F_{r-1}+m_r\kappa^F_r \] are nonzero.

We will write  $\xi\in\flagat^+$ if $m_j>\tau$, $\dt jr$, in this expression.
  \end{definition}

  If $\xi$ is a $\flaga^+$-regular element, then any flag $F$ in
  $\FFF$ is proper, and the cone $\s^+(F,\A)$ is a closed simplicial
  cone.

\begin{lemma}\label{regreg}
 \begin{enumerate}
 \item Any $\flaga^+$-regular element $\xi$ is $\A$-regular.
 \item If $\xi$ is $\A$-regular and $\flaga$-regular, then $\xi$ is
 also $\flaga^+$-regular.
 \item
 If $\xi$ is $\flaga$-regular, then for sufficiently large $s>0$,
 the element  $\xi+s\kappa_\A$ is a $\flaga^+$-regular element, moreover
$\FF=\flag^+(\A,s\kappa_\A+\xi)$.
 \end{enumerate}
\end{lemma}

\begin{remark} Note that an element $\xi$ which is
 $\flaga$-regular is not usually $\A$-regular.  For example if
 $\a^*$ has basis $e_1,e_2$ and if
 $\A:=[e_1,e_2]$ the singular element
 $\xi=e_1$ is $\flaga$-regular.
\end{remark}

\begin{proof}
First of all, let us prove by induction on the dimension of
$\ga^*$ that  for any $\xi$ in $\Cone(\A)$ there exists a proper
flag $F$ such that $\xi\in \s^+(F,\A)$. Indeed,  the point
$\kappa_\A=\kappa_r$ is in the interior of $\Cone(\A)$. Thus there
exists a face $\Face$ of $\Cone(\A)$ such that $\xi$ belongs to
the convex hull of $\kappa_\A$ and of $\Cone(\A\cap \Face)$. Let
$t\geq 0$ such that $\xi'=\xi-t\kappa_\A$ belongs to the cone
$\Cone(\A\cap \Face)$. If $F'$ is a  proper flag  for the system
$\A\cap\Face$ such that $\xi'\in \s^+(F',\A\cap\Face)$, then
$[F',\ga^*]$ is a proper flag and $\xi\in \s^+(F,\A)$.

We now prove the first statement. If $\xi$ is $\A$-singular, then
there is a hyperplane $F_{r-1}$ spanned by vectors of $\A$ such
that $\xi\in \Cone(\A\cap F_{r-1})$. We choose a flag $F'$
relative to the system $\A\cap F_{r-1}$ such that $\xi\in \s^+(F',
\A\cap F_{r-1})$. Then $F=[F',\ga^*]$ belongs to $\flaga$, $\xi$
is in $\s^+(F,\A)$, but the $\kappa_r$-coefficient of $\xi$ on
$\kappa_\A$ vanishes. This is in contradiction with  the
assumption that  $\xi$ is $\flaga^+$-regular.

To prove the second statement, observe that if $\xi$ is
$\flaga$-regular but not $\flaga^+$-regular, then $\xi$ is in a cone spanned
by the vectors $\kappa_j^F$, $\dt j{r-1}$, for some flag $F$. Clearly,
these vectors are in the cone spanned by $\A\cap F_{r-1}$. This
implies that $\xi$ is $\A$-singular.

The third item is clear.
\end{proof}

The following theorem is proved in \cite{tor1} under a stronger
assumption of regularity. However the same proof shows the
following.

\begin{theorem}[{\cite[Theorem 2.6]{tor1}}]\label{jkc}
  Let $\gc$ be any chamber of the projective sequence $\A$, and let
  $\xi$ be a vector in $\c$ which is $\flaga^+$-regular.
    Then for every $\phi\in \rata$
 \begin{equation}
   \label{jkit}
\JK_{\gc}(\phi)=\sum_{F\in \FFF} \nu(F)
 \res_{F}\phi.
 \end{equation}
\end{theorem}

Consider now a $\flaga$-regular vector $\xi$ in $\ga^*$. It follows
from Lemma \ref{regreg} that for sufficiently large real $s$, the
element $s\kappa_\A+\xi$ is $\A$-regular. We can thus associate to
every $\flaga$-regular element $\xi$ a chamber $\gc(\kappa_\A,\xi)$ as
the unique chamber containing $s\kappa_\A+\xi$ for sufficiently large
$s$. Clearly, when $\kappa_\A$ is not in $\Cone_\sing(\A)$, then the
chamber $\gc(\kappa_\A,\xi)$ is the chamber which contains
$\kappa_\A$. In general the chamber $\gc(\kappa_\A,\xi)$ is a chamber
containing $\kappa_\A$ in its closure.

Lemma \ref{regreg} and Theorem \ref{jkc} have the following corollary:
\begin{cor}\label{jk}
Let $\phi\in \rata$ and let $\xi$ be $\flaga$-regular. Then
$$\JK_{\gc(\kappa_\A,\xi)}(\phi)=\sum_{F\in \flag(\A,\xi)} \nu(F)
 \res_{F}(\phi).$$
 \end{cor}

The following key statement easily follows from Lemma \ref{regreg},
and the discussion afterwards.
\begin{prop}[{\cite[Proposition 6.1]{tor1}}]\label{vanish}
  Let $\gc$ be a chamber such that $\kappa_\A$ is in the closure of
  $\gc$, and let $\xi\in\gc$ be a $\flaga$-regular element. Then
  $\FF=\FFF$.
\end{prop}

This allows us to formulate our central result in a concise fashion as
follows:
\begin{cor}
  Let $\A=[\alpha_1,\dots,\alpha_n]$ be a projective sequence in
  $\ga^*$, and $\gc$ be a chamber of $\A$ such that
  $\kappa_\A=\sumn\alpha_i\in\bar\gc$. Then for any polynomial $Q$ on
  $\ga$ we have
\[
  \int_{V_\A(\gc)}\chi(Q)=
\int_{Z(\xi)}\frac{Q\,\dvol}{\prod_{i=1}^n\alpha_i},
\]
where $V_\A(\gc)$ is the orbifold toric variety associated to $\A$ and
$\gc$, $\chi(Q)$ is the class in $H^*(V_\A(\gc),\R)$ corresponding to
$Q$ under the Chern-Weil map, $\xi$ is a $\flagat$-regular
vector in $\gc$ and the cycle $Z(\xi)$ is defined in \eqref{defzxi}
\end{cor}

\section{The mixed mirror residue conjecture}
\label{sec:mixed}

We start with  our usual data:
$$0\to \gau\overset{\iotau}\to \ggu\overset{\piu}\to \gtu\to 0$$
an
exact sequence of real vector spaces of dimensions $\ru, \nun, \du$
respectively; each vector space is endowed with an integral structure,
and the sequence is exact over the integers. The vector space $\ggu$
is endowed with a basis $\{\omega_i, 1\leq i\leq \nun\}$, so that
$$\ggu=\oplus_{i=1}^{\nun} \R
\omega_i,\quad\text\quad\ggu_\Z=\oplus_{i=1}^{\nun} \Z \omega_i.$$

Let $\Bu=[\betau_1,\betau_2,\ldots, \betau_{\nun}]$ be the sequence of
vectors in $\gtu_{\Z}$ which are images of the basis vectors under
$\piu$: $\betau_i=\piu(\omega_i)$. We assume that $0$ is in the
interior of the convex hull $\Pi^{\Bu}$ of the elements of $\Bu$, and
that $\Bu$ generates $\gt_\Z$ over $\Z$. We denote by
$\Au=[\alphau_1,\alphau_2,\ldots, \alphau_{\nun}]$ the Gale dual
sequence in $\gau^*$, i.e. the sequence obtained as the restrictions
to $\gau$ of the coordinate functions $\omega^i\in\ggu^*$ ($1\leq
i\leq \nun$). Then $\Au$ generates $\gau^*$.  We introduce the
notation $\kappau$ for the vector
$\kappa_{\Au}=\sum_{i=1}^{\nun}\alphau_i$.

Following \cite{BMmixed},  consider a decomposition of the set
$\{1,2,\ldots,\nun\}$   into $l$ disjoint sets:
\begin{equation}
  \label{partD}
\{1,2,\ldots,\nun\}=\bigcup_{k=1}^l D_k,
\end{equation}
and define $\theta_k=\sum_{i\in D_k}\alphau_i$. Note that
$\kappau=\sum_{k=1}^l\theta_k$.

 Given this data,
we construct an example of the setup described at the beginning of
the paper.

We set  $V=\gtu\oplus \R^{l}$; let
$[\gamma_1,\ldots,\gamma_l]$ be the canonical basis of $\R^l$, so that
$$V=\gtu\oplus \bigoplus_{k=1}^l \R\gamma_k.$$
Then $V$ is a $d+1=\du+l$ dimensional vector space endowed with a
lattice.

We choose the grading vector $g\in V^*$ to be the linear functional
which vanishes on $\gtu$ and takes the value 1 on the vectors
$\gamma_k$, $k=1,\dots,l$. The set $\M$ is constructed as follows: for
each $\betau_i\in \Bu$, let $\mu_i=\betau_i+\gamma_{k}\in V$, where
$k$ is the integer satisfying $i\in D_k$, and set
\[ \mu_i=\gamma_{i-\nun}, \quad i=\nun+1,\dots,\nun+l.\] Thus our sequence
$\M=[\mu_1,\dots,\mu_n]$ contains $n=\nun+l$ elements, and further
defines an exact sequence
\[ 0\to W\to\ggu\oplus\R^l\to\gtu\oplus\R^l\to 0,\]
where the map from $\ggu\oplus \R^l$ to $\gtu\oplus\R^l$ sends
$\omega_i$ with $i\leq \nun$ to $\mu_i=\betau_i+\gamma_k$ and
$\gamma_k$ to $\gamma_k$.  Comparing this with the diagram
\eqref{diagram}, we see that we have $\gg=\ggu\oplus\R^l$, and
$W=\gau$.

We can identify the space $W\subset\ggu\oplus\R^l$ with $\gau$ via the
map $u\mapsto u-\sum_{k=1}^l\langle \theta_k,u\rangle \gamma_k$.
Indeed, the image of the element
$u=\sum_{i=1}^{\nun}\alphau_i(u)\omega_i\in \ggu\oplus\R^l$ in
$\gtu\oplus\R^l$ is
$\sum_{i=1}^{\nun}\alphau_i(u)(\betau_i+\gamma_k)=\sum_{k=1}^l\langle
\theta_k,u\rangle \gamma_k$.

We observe that the cone of $C=C(\M)$ defined in \S\ref{sec:notation}
is exactly the construction of the Cayley cone given in
\cite[\S3]{BMmixed}.  Denote by $\CH_k$ the convex hull of $\{0\}$ and
the elements $\beta_i$ with $i\in D_k$. Then the base of the cone
$C=C(\M)$ is the convex hull of the union
$$\CH_1\times \{\gamma_1\}\cup \cdots
\cup \CH_l\times \{\gamma_l\}\subset\gg,$$ and our cone $C(\M)$ is the cone
$C(\CH_1,\CH_2,\ldots, \CH_r)$ in \cite[Definition 3.1]{BMmixed}.

Implementing the construction in \S\ref{sec:notation} further, we see
that in our present setup there is a natural candidate for the vector
$\gamma\in V$:
\begin{equation}
  \label{defgam}
\gamm=\sum_{k=1}^l\gamma_k.
\end{equation}
This vector indeed lies in the interior of the cone $C$ generated by
$\M$; it has degree $l$.

Now we define the vector space $\gt=V/\R \gamm$, which is of dimension
$d=\du+l-1$, and we consider the map $\pi:\gg\to \gt$,
which sends ${\omega_i}$ to $\mu_i$ modulo $\gamm$. The images of the
basis vectors will again be denoted by $\beta$s, thus we have
\[ \B =[\beta_i=\pi(\mu_i),\;i=1,\dots, n].
\]

Note that $\pi(\gamm)=0$ and that the space $\ga$, the kernel of
$\pi$, may be canonically identified with $\gau\oplus\R\gamm$. This
decomposition holds over the integers as well.

Now denote by $\deltau$ the linear functional on $\ga$ which is zero
on $\gau$ and satisfies $\lr\deltau\gamma=1$. Then we have a canonical
decomposition $\ga^*=\gau^*\oplus\R\deltau$. Again, we compute the Gale
duals, and we obtain
\begin{lemma}\label{desca}
 The Gale dual sequence $\A$ of the sequence $\B$ is
\[ \alphau_i,\, i=1,\dots,\nun;\;
\alpha_i=\deltau-\theta_{i-\nun},\; i=\nun+1,\dots,n.\] Thus
$$\A=[\alphau_1,\alphau_2,\ldots,\alphau_{\nun}, \deltau-\theta_1,
\deltau-\theta_2,\ldots, \deltau-\theta_l].$$

\end{lemma}
Note that since $\sum_{k=1}^l \theta_k=\sum_{i=1}^{n}\alphau_i$,
we have $\kappa_\A=l\deltau$.

\bigskip

Consider the partition polytope $\Pi_{\kappau}\subset \ggu^*$ and the
point $b=\sum_{i=1}^n \omega^i\in \Pi_{\kappau}$.  Define the polytope
$\PPu$ as $\Pi_{\kappau}-b$. An easy computation shows that the dual
polytope of $\PPu$ is the convex hull of the set $\Bu$.  Consider the
partition polytopes $\Pi_{\theta_k}$, and define
$\PPu_k=\Pi_{\theta_k}-b_k\subset \gtu^*$, where $b_k=\sum_{i\in
  D_k}\omega^i\in\Pi_{\theta_k}$.

The following proposition follows directly from Proposition
\ref{minkowski}.

\begin{prop}

 The  polytope $\PPu$ is the Minkowski sum
of the polytopes $\PPu_k$ if  and only if the elements $\theta_k$
are in the closure of the same chamber $\gc\subset \gau^*$.
\end{prop}

\begin{definition}
Let $\gc$ be a chamber of the $\Cone(\Au)$. We will say that the
projective sequence $\Au$ is $\gc$-partitioned into $l$ parts if the
chamber $\gc$ contains in its closure all the elements
$\theta_k=\sum_{i\in D_k}\alphau_i$, $k=1,\ldots,l$.
\end{definition}
Note that in this case the chamber $\gc$ also contains in its closure
the element $\kappau=\kappa_{\Au}=\sum_{k=1}^l\theta_k$.

Consider (cf. \cite[Definition 4.2]{BMmixed}) the dual family
$\PPo_1,\ldots,\PPo_l$  of the family of polytopes  $\PPu_k$ defined by
the equations
$$ \PPo_k=\{y\in \gtu \, ; \langle x, y \rangle \geq \-\delta_{ks},
x\in \PPu_s, s=1,\ldots, l\}.$$

\begin{prop}( \cite[Definition 4.2]{BMmixed})
Assume $\Au$ is $\gc$-partitioned into $l$ parts. Then the
polytope $\PPo_k$ is the convex hull of the origin $\{0\}$ and the
vectors $\betau_i$ with $i\in D_k$.
\end{prop}

\begin{proof}
  First we verify that the points $\betau_i$ with $i\in D_k$ are in
  $\PPo_k$. Let $c\in \PPu_s\subset \gtu^*$. Then $c$ may be written
  as $\sum_{i=1}^{\nun}c_i \omega^i -\sum_{i\in D_k} \omega^i$, with
  $c_i\geq 0$. Since $\beta_i=\pi(\omega_i)$, the value of $c$ on
  $\beta_i$ is $\langle c,\omega_i\rangle$; this, in turn, equals
  $c_i-1$ if $i\in D_k$, and $c_i$ if $i\notin D_k$. Thus $\beta_i$
  with $i\in D_k$ satisfies the inequalities defining $\PPo_k$.

  Conversely, let $y\in \PPo_k$, with $y\neq 0$, and choose a cone of
  the fan in $\gtu$ determined by $\gc$, which contains $y$. This
  means that there is a $\du$-element subset
  $\nu\subset\{1,\dots,\nun\}$ such that $y=\sum_{i\in \nu}t_i
  \beta_i$ with $t_i\geq0$, and $\gc\subset \sum_{i\notin\nu}
  \R^+\alpha_i$. Since we have a $\gc$-partition, this implies that we
  can write $\theta_s=\sum_{j\notin \nu} x_s^j \alpha_j$, with
  $x_s^j\geq 0$; this means that the point $\rho_s=\sum_{j\notin \nu}
  x_s^j \omega^j -\sum_{i\in D_s}\omega^i$ is in $\PPu_s$ for $1\leq
  s\leq l$.  By the definition of $\PPo_k$, we have
  $\lr{\rho_s}y\geq0$ for $s\neq k$.  Since $\lr{\rho_s}y\geq0=\sum_{i\in
    D_s\cap \nu}t_a$, it follows that $\nu$ is contained in $D_k$.
  We also have $\lr{\rho_k}y\geq-1$, hence
  $-\sum_{i\in \nu}t_i\geq -1$, thus $y=\sum_{a\in \nu} t_i \beta_i
  +(1-\sum_{i\in \nu}t_i)0$ is in the convex hull of the elements
  $\{\beta_i,\,i\in D_k\}$ and $0$.
\end{proof}

The following proposition will be crucial in the proof of the
mixed toric residue mirror conjecture.

\begin{prop}\label{descflag}
  Assume $\Au$ is $\gc$-partitioned into $l$ parts, i.e. $\theta_k\in
  \bar{\gc}$, $k=1,\ldots,l$. Let $\xiu\in\gc\subset\gau^*$ and
  $\xi\in\ga^*$ be such that $\xi-\xiu=b\deltau$ for some $b\in \R$.
  Then the map $\Fu \mapsto (\Fu,\ga^*)$ from $\flag(\Au)$ to
  $\flag(\A)$ induces a bijection between $\flag(\Au,\xiu)$ and $
  \flag(\A,{\xi})$.
\end{prop}

\begin{proof}
  The injectivity of the correspondence $\Fu \mapsto (\Fu,\ga^*)$ is
  obvious, thus we only need to show surjectivity. Let
$$F = [F_0=\{0\}\subset F_1\subset F_2\subset  \dots \subset
F_{r-1} \subset F_r \subset
F_{r+1}=\gau^*\oplus\R\deltau]\in\flag(\A,\xi).$$

By definition of $\flag(\A,\xi)$, we can write
$\xi=\sum_{j=1}^{r+1} c_j
  \kappa_j^F $ with $c_j\geq 0$ for $j\leq r$.
Using the description of $\A$ in Lemma \ref{desca} and projecting
this
  equality onto $\gau^*$ along $\R\deltau$ we obtain
\begin{equation}\label{singular}
\xiu+\sum_{j=1}^rc_j\sum_{\{k;\,t-\theta_k\in F_j\}}\theta_k=
\sum_{j=1}^rc_j\sum_{\{i;\,\alphau_i\in F_j\}}\alphau_i.
\end{equation}
Now assume that $F_r\neq\gau^*$. Then the elements $\alphau_i$ on
the right hand side do not span $\gau^*$, and, therefore, the
right hand side lies in $\Cone_{{\rm sing}}(\Au)$.  On the other
hand, according to our assumptions on the $\theta_k$s and $\xiu$,
the left hand side of \eqref{singular} is in $\gc$. This is a
contradiction, since a vector is in $\Cone_{{\rm sing}}(\Au)$
exactly when it is not in a chamber. Thus we have shown that any
flag
  $F\in \flag(\A,{\xi})$ is of the form $ F=(\Fu,\ga^*)$, where $\Fu$
  is an element of $\flag(\Au,\xiu)$.
\end{proof}
\begin{cor}\label{iffreg}
  With the assumptions and notations of the Proposition, the vector
  $\xiu\in\gc\subset\gau^*$ is $\flagau^+$-regular
  ($\flagau_\tau^+$-regular) if and only if the vector $\xi+bt\in\ga^*$
  is $\flaga$-regular ($\flagat$-regular).
\end{cor}

Now let $z=(z_1,z_2,\ldots,z_{\nun})$ be a vector in $\C^{\nun}$, and
introduce the vector $\tz=(z_1,z_2,\ldots,
z_{\nun},1,1,\ldots,1)\in\C^n$.

Then, as in \S\ref{sec:torres}, we consider
$$f_\tz=\sum_{k=1}^l e_{\gamma_k}+\sum_{i=1}^{\nun} z_i
e_{\mu_i}\in S^1(C).$$
This function may be written as
$$f_\tz=\sum_{k=1}^l e_{\gamma_k}(1+\sum_{i\in D_k} z_i e_{\betau_i}).$$
Note that after the substitution $z_i=-a_i$, our function $f_\tz$ is the
Cayley polynomial $F$ associated to the functions $f_k$ and the
sets $A_k=\{\beta_i; i\in D_k\}$ in \cite[\S4]{BM}.

As in \S\ref{sec:torres}, set $f_{\tz,0}=f_\tz$ and consider the
sequence $\{f_{\tz,1},\dots, f_{\tz,d}\}$ of derivatives of the
function $f_\tz$, with respect to a basis $\{a_1,a_2,\ldots, a_d\}$ of
$\gt^*$.  Comparing this setup to \S\ref{sec:torres} on toric residue,
we see that the function $f_\tz$ is the function associated to the
sequence $\M$ and the {\em non-generic} parameter $\tz$.  However, we
have the following lemma.

\begin{lemma}\label{condlemma}
  For generic $z\in\C^{\nun}$, the sequence $\M$, the vector $\gamma$
  given in \eqref{defgam} and the parameter $\tz$ satisfy conditions
  \eqref{cond1} and \eqref{cond2}. This means that the sections
  $\{f_{\tz,0},f_{\tz,1},\ldots, f_{\tz,d}\}$ have no common zeroes in
  $\Toric(C)$.  Furthermore, the common zeroes of the sections
  $\{f_{\tz,1},\ldots, f_{\tz,d}\}$ are in the torus
  $\Torus(H^*)\subset\Toric(C)$.
\end{lemma}
\begin{proof}
  This is a consequence of the results of the Appendix. Indeed
  consider the subset $\sigma\subset \{1,2,\ldots,n\}$ consisting of
  $\nun+1, \ldots,n=\nun+l$. Then $(\C^*)^{\nun}$ may be identified with the
  set $Z_\sigma$ of Definition \ref{Zsigma} by associating to the
  vector $z=(z_1,\ldots, z_{\nun})$ the vector $\tilde z=(z_1,\ldots,
  z_{\nun},1,\ldots,1)$. The corresponding $\mu_i$ with $i\in \sigma$
  are the set of $l$ vectors $\gamma_k$, $k=1,\ldots,l$. These are
  linearly independent by definition, and thus the statement of the
  Lemma follows from Propositions
  \ref{denseopen} and \ref{possible} of the Appendix.
\end{proof}

Thus for generic $z$ in $\C^{\nun}$ the conditions \eqref{cond1}
and \eqref{cond2} are satisfied and the toric residue with respect
to $\bftz$ is well defined. Moreover, the conditions of Proposition
\ref{local} are also satisfied.

Now we are ready to formulate our main result: a generalization of the
mixed mirror residue conjecture of Batyrev and Materov.
Fix a polynomial $P$ of degree $\du$ in $\nun$ variables. Then
\begin{equation}
  \label{thep}
e_{\gamm}P(z_1 e_{\mu_1},\ldots, z_{\nun}e_{\mu_{\nun}})
  \end{equation}
  is an element of $I^{d+1}(C)$.  The toric residue of this element
  with respect to the sequence $\bftz$ is well-defined; it is a rational
  function of $z$.  The conjecture expresses this toric residue as a
  sum of intersection numbers on a sequence of Morrison-Plesser moduli
  spaces.

  Assume that we have a $\gc$-partition of  $\Au$, i.e. that all the
  elements $\theta_k\in \gau^*$, $k=1,\ldots,l$ belong to the closure
  of the same chamber $\gc$ of the projective sequence $\Au$.
  Consider the polar cone $\gcd$
\[ \gcd = \{\lambda\in\gau;\; \lr\alpha\lambda\geq0\text{ for all
}\alpha\in\gcu\}.\] If $\lambda\in \gcd$, then the numbers
$\langle\theta_k,\lambda\rangle $ are non-negative integers. In
particular,
$$
\frac{P(\alphau_1,\dots,\alphau_{n}) \prod_{k=1}^l
  \theta_k^{\langle\theta_k,\lambda\rangle }}
{\prod_{i=1}^{n}\alpha_i^{\langle \alpha_i,\lambda\rangle +1}}$$
is an element of degree $-r$ of $\C_{\Au}[\gau]$, and we can take its
Jeffrey-Kirwan residue with respect to the chamber $\gc$.

\begin{theorem}\label{bat}
  Let $\gc$ be a chamber of the projective sequence $\Au$, and assume
  that $\Au$ is $\gc$-partitioned into $l$ parts.  Let $C=C(\M)$ be
  the corresponding Cayley cone and let $f_\tz=\sum_{k=1}^l
  e_{\gamma_k}(1+\sum_{i=1}^{\nun} z_i e_{\betau_i})$ be the Cayley
  polynomial associated to the partition. Denote by $\bftz$ the
  associated sequence of partial derivatives of $f_tz$.

  Choose a vector $z\in \C^{\nun}$ such that $\psi(z)=\psi_{\Au}(z)$
  is a $\flagaut$-regular element of $\gc$ for some sufficiently large
  $\tau$.  Let $P$ be a polynomial of degree $\du$ in $\nun$
  variables.  Then the series
  \begin{equation}
    \label{theseries}
\sum_{\lambda\in \gcd}
  \JK_{\gc}\left(\frac{P(\alphau_1,\dots,\alphau_{\nun})
  \prod_{k=1}^l\theta_k^{\langle\theta_k,\lambda\rangle }}
  {\prod_{i=1}^{\nun}\alphau_i^{\langle \alphau_i,\lambda\rangle +1}}
  \right)z^{\lambda}
  \end{equation}
  converges absolutely, and sums to the toric
  residue
\begin{equation}\label{thetoricres}
\tres_{\bftz}e_{\gamm}P(z_1e_{\mu_1},\ldots,z_{\nun}e_{\mu_{\nun}})
\end{equation}
\end{theorem}

\begin{remark}

\begin{itemize}
\item In \cite{BMmixed}, Batyrev and Materov requires that all
  polytopes $\PPu_k$ have integral vertices.  We do not require this
  hypothesis.
\item
If the polynomial $P$ is supposed to be multi-homogeneous of
degree $\du_k$ with respect to the variables $x_i$ with $i\in D_k$
with $\sum_{k=1}^l \du_k=\du$, then
$$P(z_1e_{\mu_1},\ldots,z_{\nun}e_{\mu_{\nun}})=
e_{\du_1\gamma_1}\cdots e_{\du_k \gamma_k}P(z_1
e_{\betau_1},\ldots, z_{\nun}e_{\betau_{\nun}})
$$
and $e_{\du_1\gamma_1}\cdots e_{\du_k \gamma_k}P(z_1
e_{\beta_1},\ldots, z_{\nun}e_{\beta_{\nun}})$ is an element of
$S^{\du}(C)$. In the situation considered in \cite{BMmixed}, the
Cayley cone $C$ is a reflexive Gorenstein cone and $I^{\du+l}(C)$
is isomorphic to $S^{\du}(C)$ by multiplication by $e_\gamma$.

\item
As shown in \cite{tor1} the number $I(P,\beta)$ of \cite{BMmixed}
is equal to
$$\JK_{\gc}\left(\frac{P(\alphau_1,\dots,\alphau_{\nun})
\prod_{k=1}^l\theta_k^{\langle\theta_k,\lambda\rangle }}
  {\prod_{i=1}^{\nun}\alphau_i^{\langle \alphau_i,\lambda\rangle+1}}\right)$$
\end{itemize}
\end{remark}
\begin{proof}
  We showed in Lemma \ref{condlemma} that the conditions of the
  Proposition \ref{local} hold in the case of the toric residue
  \eqref{thetoricres}. As a result we can apply the formula
  \eqref{torlocal}.
  
  As explained in section \S\ref{sec:gale}, we can parameterize the
  set $\tilCrit(z)$ by elements
\[ (u,t)\in  O(z,\A)\subset \gau_\C\oplus \C\gamm.\]

Note that the character $e_\gamma$ is equal to $1$ on $O(z,\A)$, thus
for $(u,t)\in O(z,\A)$, we have
$\prod_{k=1}^l(\deltau-\theta_k(u))=1.$

 Recall that $\kappa(\A)=l\deltau$.  Then the toric residue
 \eqref{thetoricres} may be written as
 $$l\sum_{(u,t)\in O(z,\A)}\frac{P(\alphau_1(u),\ldots,
 \alphau_{\nun}(u))}{l\deltau G(u,t)},$$
where $G(u,t)$ is given by \eqref{G}.

Recall that we associated the vector $\tz=(z_1,\ldots,
z_{\nun},1,\ldots,1)\in\C^n$ to $z\in\C^{\nun}$. Then
$\psi_{\A}(\tilde z)=(\psi_{\Au}(z),0)\in \ga^*=\gau^*\oplus \C
t$. Recall from Corollary \ref{iffreg} that the vector
$\psi_{\A}(\tilde z)\in \ga^*$ is $\flagat$-regular if and only if
$\psi_{\Au}(z)\in \gau^*$ is ${\mathcal{FL}}^+(\Au)$-
$\tau$-regular.  We choose $\tau$ large enough in order that
Theorem \ref{FORMULA} holds for $\A$.

Then using Lemma \ref{desca} and Theorem \ref{FORMULA} we can conclude
that the toric residue \eqref{thetoricres} is equal to
$$\sum_{F\in \flag(\A,\xi)}\nu(F) \sum_{\Lambda\in
  \ga_\Z}\res_{ F}\left(
  \frac{K(u,\deltau)}{p_\Lambda(u,\deltau)}\frac{1}{\prod_{i=1}^{\nun}
    \alphau_i(u)\prod_{k=1}^l(\deltau-\theta_k(u))}\right) \tilde
z^{\Lambda}$$
with
$$K(u,\deltau)=\frac{P(\alphau_1(u),\ldots,
\alphau_{\nun}(u))}{\deltau}.$$

We write $\Lambda\in \ga_\Z$ as $\Lambda=\lambda+m \gamm$ with
$\lambda\in \ga_\Z$ and $m\in \Z$. Then we have
$$p_\Lambda(u,\deltau)= p_\lambda(u)\prod_{k=1}^l(t-\theta_k(u))^{m-\langle
\theta_k,\lambda\rangle },$$ while $(\tilde z)^{\Lambda}$ is
simply $z^{\lambda}$. To simplify our notation, we introduce
 $$P_{\Au}(u)=P(\alphau_1(u),\ldots,
\alphau_{\nun}(u)).$$ This is a polynomial on $\gau_{\C}$ of degree
$d$, but it can be considered as a polynomial on $\ga_{\C}$ as well.

Thus, after some simplifications, changing $m$ to $m-1$ and using
the absolute convergence of the series, we obtain that the toric
  residue \eqref{thetoricres} is equal to
$$\sum_{\lambda\in \ga_\Z,m\in
  \Z}\sum_{ F\in \flag(\A,\xi)}\nu(F) \res_{ F}\left(
  \frac{P_{\Au}(u)}{\deltau
    p_\lambda(u)\prod_{k=1}^l(\deltau-\theta_k(u))^{m-\langle
      \theta_k,\lambda\rangle }}\frac{1}{\prod_{i=1}^{\nun}
    \alphau_i(u)}\right) z^{\lambda}.$$

According to Proposition \ref{descflag}, for $\xiu\in\gc$ and
$\xi=\xiu+b\deltau$,  for any $b\in \R$, in particular for $b=0$,
we have
\begin{equation} \label{flagspm}
 \flag(\A,\xi)=\{(\Fu,\gau^*);\;\Fu\in\flag(\Au,\xiu)\}
\end{equation}

Let us take $\Lambda=\lambda+m\gamma\in \ga_\Z$, fix
$F\in\flag(\A,\xi)$, and denote by $a(F,\Lambda)$ the
corresponding term in the above sum, omitting the sign $\nu(F)$.
Thus
$$a(F,\Lambda)=\res_{ F} \left(\frac{P_{\Au}(u)}{\deltau p_\lambda(u)
    \prod_{k=1}^l(\deltau-\theta_k(u))^{m-\langle
\theta_k,\lambda\rangle}}\frac{1}{\prod_{i=1}^{\nun}\alphau_i(u)}\right).$$
Then translating the
relation \eqref{flagspm} into the language of iterated residues, we
obtain
\begin{equation} \label{afl}
a(F,\Lambda)=\res_{\deltau=0}\frac{d\deltau}{\deltau}\res_{\Fu}
\left(\prod_{k=1}^l (\deltau-\theta_k(u))^{\langle \theta_k,\lambda\rangle
    -m}\frac{P_{\Au}(u)}
  {p_\lambda(u)\prod_{i=1}^{\nun}\alphau_i(u)}\right).
\end{equation}

\begin{lemma}\label{mzero}
  Let $\Lambda=\lambda+m\gamma\in \ga_\Z$. Then $a(F,\Lambda)=0$ unless
  $m=0$.
\end{lemma}
\begin{proof}
  The proof will follow from degree considerations. Since $t=0$ is the
  last residue, in order to compute \eqref{afl}, we need to write
 \begin{equation}\label{expand}
\prod_{k=1}^l (\deltau-\theta_k(u))^{\langle
    \theta_k,\lambda\rangle -m}=\deltau^{-lm+\langle
    \kappau,\lambda\rangle } \prod_{k=1}^l
  \left(1-\frac{\theta_k(u)}{\deltau}\right)^{\langle \theta_k,\lambda\rangle
    -m}
\end{equation}
  and expand $(1-\theta_k(u)/\deltau)^{\langle
    \theta_k,\lambda\rangle -m}$ in a power series involving positive
  powers of $(\theta_k(u)/\deltau)$. We will track $\deg_\deltau$, the
  $\deltau$-degree, and $\deg_{\Au}$, the degree of the rational
  function on $\gau$. Clearly, after the expansion of the expression in
  parentheses in \eqref{afl} we will need to retain the terms with

\begin{equation}\label{degrees}
\deg_\deltau=0\quad \text{and}\quad \deg_{\Au}=-\ru.
\end{equation}
Now note that $\deg_{\Au} P_{\Au}=\du$,
  $\deg_{\Au}p_\lambda=\lr\kappau\lambda$, and the difference of
  degrees $\deg_\deltau-\deg_{\Au}$ of the expression \eqref{expand}
  is equal to $\lr\kappau\lambda-lm$. Hence we can conclude that the
  difference of degrees
\[ [\deg_\deltau-\deg_{\Au}]\left(\prod_{k=1}^l (\deltau-\theta_k(u))^{\langle \theta_k,\lambda\rangle
    -m}\frac{P_{\Au}(u)}
  {p_\lambda(u)\prod_{i=1}^{\nun}\alphau_i(u)}\right)
\] of the expression in \eqref{afl} is equal to $\ru-lm$.
Comparing this to \eqref{degrees} we see that the only possibility is $m=0$.
\end{proof}

Now we are ready to finish the proof of Theorem \ref{bat}. Taking
into account Lemma \ref{mzero}, and switching back the order of
summation, we have
\begin{multline*}
\tres_{\bftz}e_{\gamm}P(z_1e_{\mu_1},\ldots,z_{\nun}e_{\mu_{\nun}}) =
\sum_{\lambda\in
\gau_\Z}\\\sum_{\Fu\in \flag(\Au,\xiu)}\nu(\Fu)
\res_{\deltau=0}\frac{d\deltau}{\deltau}\res_{\Fu}
\left(\prod_{k=1}^l (\deltau-\theta_k(u))^{\langle \theta_k,\lambda\rangle
    }\frac{P_{\Au}(u)}
  {p_\lambda(u)\prod_{i=1}^{\nun}\alphau_i(u)}\right) z^{\lambda}.
\end{multline*}
We observe that since $\xiu$ is in chamber $\gc$ containing
$\kappau=\kappa(\Au)$ in its closure, according to Proposition
\ref{vanish} and Theorem \ref{jkc}, we can replace the sum of iterated
residues by a Jeffrey-Kirwan residue, and we arrive at the expression
\[ \sum_{\lambda\in
\gau_\Z}
\res_{\deltau=0}\frac{d\deltau}{\deltau}
\JK_\gc
\left(\prod_{k=1}^l (\deltau-\theta_k(u))^{\langle \theta_k,\lambda\rangle
    }\frac{P_{\Au}(u)}
  {p_\lambda(u)\prod_{i=1}^{\nun}\alphau_i(u)}\right) z^{\lambda}.
\]
for the toric residue.

Next, using Lemma \ref{null}, we can restrict the summation on the
right hand side to the integral vectors $\lambda$ in $\gcd$. Now we
use our assumption that the vectors $\theta_k$ are in $\bar\gc$, which
implies that $\lr{\theta_k}\lambda\geq0$ for every $\lambda\in\gcd$ and
$k=1,\dots,l$. This means that all the exponents in the product
\[ \prod_{k=1}^l (\deltau-\theta_k(u))^{\langle
  \theta_k,\lambda\rangle}
\]
are nonnegative, which immediately implies that only the term
\[ \prod_{k=1}^l \theta_k(u)^{\langle
  \theta_k,\lambda\rangle},
\]
can contribute to the residue with respect to $\deltau$. Thus we are
left with the equality
\begin{multline*}
\tres_{\bftz}e_{\gamm}P(z_1e_{\mu_1},\ldots,z_{\nun}e_{\mu_{\nun}})=\\
\sum_{\lambda\in \gau_\Z\cap\gcd} \JK_\gc
\left(\prod_{k=1}^l \theta_k(u)^{\langle \theta_k,\lambda\rangle
    }\frac{P_{\Au}(u)}
  {p_\lambda(u)\prod_{i=1}^{\nun}\alphau_i(u)}\right) z^{\lambda}.
 \end{multline*}
This completes the proof.

\end{proof}

\section{Appendix}
\label{sec:appendix}

We start again with the setting of \S\ref{sec:localform}. Let
$z=(z_1,z_2,\ldots, z_n)$ be a vector in $\C^n$, and consider the
element
$$f_z=\sum_{i=1}^n z_i e_{\mu_i}\in S^1(C).$$
 For each $a\in V^*$, we consider the derivative
$$f_{z,a}=\sum_{i=1}^n z_i \langle a,\mu_i\rangle e_{\mu_i}$$
of the function $f_z$ in the direction of the vector $a\in V^*$.

Choose a basis $[a_0,a_1,\ldots, a_d]$ of $V^*$, where $a_0=g$, the
grading vector.  Our goal is to describe a suitable set $Z$ of values
of the vector parameter $z$ for which the sections $f_{z,j}=f_{z,a_j}$ do
satisfy the conditions \eqref{cond1} and \eqref{cond2}. The first
condition is
\begin{eqnarray}\label{cond11}
 \{x\in\Toric(C);\; f_{z,a_j}(x)=0, \,j=0,\dots,d\}=\emptyset
\end{eqnarray}
when $z\in Z$. Note that this condition does not depend of the choice
of the basis $[g,a_1,\ldots,a_d]$.  In this appendix, we recall an
argument we learned from Alicia Dickenstein to determine such a
suitable set $Z$ of values of the vector parameter $z$.

Recall the notations of \S\ref{sec:gale}. Let  $\ggy=\R^n.$ We denote by
$\omega_i$ the canonical basis of $\R^n$ so that $$\ggy=
\oplus_{i=1}^n \R \omega_i.$$  We introduce the corresponding
lattice
$$\ggy_\Z=\oplus_{i=1}^n \Z \omega_i.$$
The dual basis $\omega^i\in
\ggy^*$ is the set of coordinates on $\ggy$.  Consider the surjection
$\ggy\to V$ sending $\omega_i$ to $\mu_i$. Our grading vector $g\in
V^*$ takes value $1$ on each $\mu_i$; let $\gr\subset \ggy$ be the
kernel of this map. We then obtain an exact sequence
$$0\to  \gr\to \ggy\to V\to 0.$$
We denote by $\gr_\Z$ the intersection of the lattice $\ggy_\Z$
with $\gr$. We denote by $\nu_i$ the restriction to $\gr$ of the
coordinate $\omega^i$.  An element $w\in \gr$ is written in the
basis $\omega_i$ as  $w=\sum_{i=1}^n\nu_i(w)\omega_i$. By
definition, for all $w\in \gr$ and $a\in V^*$, we have the
relation
\begin{equation}\label{R1}
\sum_{i=1}^n \nu_i(w)\mu_i(a)=0.
\end{equation}

Note that if  $m=\sum_{i=1}^n m_i\omega_i\in \gr$ we have
$\sum_{i=1}^nm_i=0$, as follows from the relation $\sum_{i=1}^n
m_i \mu_i(g) =0$.

\bigskip
We start our analysis of the common set of zeroes of the functions
$f_{z,a_j}$.

Let  $x\in \Spec(S(C))$ be an element of the affine variety
$\Aff(C)$.  Recall the notation $x_i$ for $x(e_{\mu_i})$.

\begin{lemma}
 If $x\in \Aff(C)$ is
 a common zero of the functions
$f_{z,a_j},j=0,\ldots,d$, we have the relation
$$\sum_{i=1}^n z_i x_i
\omega_i\in \gr_\C.$$
\end{lemma}

\begin{proof}  Considering the value of $x$ on $f_{z,a_j}\in S^1(C)$,
 we obtain the set of equations $\sum_{i=1}^n z_i \langle
a_j,\mu_i\rangle x_i=0$ where  $a_j$ runs through a basis of
$V^*$. Thus $\sum_{i=1}^n z_i x_i \mu_i=0$, which exactly says
that $\sum_{i=1}^n z_i x_i \omega_i$ belongs to $\gr_\C$.

\end{proof}

Now we analyze conditions on the vector parameter $z$ which guarantee
that the equation $\sum_{i=1}^n z_i x_i \omega_i\in \gr_\C$
together with the fact that $x$ belongs to $\Aff(C)$ implies that
$x=0$. We  need to recall the equations of $\Aff(C)$.

For every element $m=\sum_{i=1}^n m_i \omega_i\in \gr_{\Z}$, we
have $\sum_{i=1}^n m_i\mu_i=0$, which, in turn, implies the
following relation in $S(C)$:
\begin{equation}\label{R2}
\prod_{\{i;\,m_i>0\}}(e_{\mu_i})^{m_i}=
\prod_{{\{i;\,m_i<0\}}}(e_{\mu_i})^{-m_i}.
\end{equation}

This gives us the set of binomial homogeneous equations on $x\in
\Aff(C)$:

\begin{equation}\label{R3}
\prod_{\{i;\,m_i>0\}}x_i^{m_i}= \prod_{{\{i;\,m_i<0\}}}x_i^{-m_i}
\end{equation}

We identify $\ggy_\C$ with $\C^n$ with the help of the basis
$\omega_i$. For every $m=\sum_{i=1}^n m_i \omega_i\in \ggy_\Z$, we
introduce the rational function $p_{m}(z)=\prod_{i=1}^n z_i^{m_i}$
on $\ggy_\C$, which is a Laurent monomial of homogeneous degree
$\sum_{i=1}^n m_i$. It is well defined on $(\C^*)^n$. Note that
the rational functions $p_m(z)$ have the multiplicative property
$$p_{m_1} p_{m_2}=p_{m_1+m_2}$$ for $m_1,m_2\in \ggy_\Z$.

Equation \eqref{R3} shows that  if $x\in \Aff(C)$ is such that $x_i$
is non zero for all $i$, then for all $m\in \gr_{\Z}$, the number
$p_m(x_1,x_2,\ldots, x_n)=\prod_{i=1}^n x_i^{m_i}$ is equal to
$1$. Our strategy is to show that there exists a non-empty finite
set $S$ in $\gr_\Z$ and constants $c_m\neq 0$ such that the
relation $\sum_{m\in S} c_m p_m=0$ holds identically on $\gr_\C$.
If we evaluate $\sum_{m\in S} c_m p_m $ on the element
$\sum_{i=1}^n z_i x_i \omega_i$ of $\gr_\C$, we obtain the
relation $\sum_{m\in S} c_m p_m(z_1,z_2,\ldots, z_n)=0$. In
particular, if our vector $z\in \C^n$ is chosen outside the
hypersurface $\sum_{m\in S} c_m p_m=0$, the equation $\sum_{i=1}^n
z_i x_i\omega_i\in \gr_\C$ cannot hold. We will need to refine
this argument to obtain also a contradiction, when some of the
$x_i$ vanishes.

For $x\in \Aff(C)$, let us analyze the set $I(x)\subset
\{1,2,\ldots,n\}$ of elements $i$ with $x_i\neq 0$.

\begin{lemma}\label{face}
Let $x\in \Aff(C)$. If $x\neq 0$, then there exists a face $\Face$ of
the cone $C$ such that $x_i\neq 0$ if and only if $\mu_i\in \Face$.
\end{lemma}

\begin{proof}
If $x\neq 0$, the set $I:=I(x)$ is not empty. Consider the cone
$C(I)$ generated by the elements $\{\mu_i;\, i\in I\}$. Let us
prove that this cone is a face of the polyhedral cone $C$.
Consider a minimal face $\Face$ among the faces containing $C(I)$.
Then $C(I)$ contains an interior point $\mu$ of $\Face$. We can
write $N\mu=\sum_{i\in I} n_i \mu_i$ with  $N$ a positive integer
and $n_i\geq 0$.  The  relation $e_{N\mu}=\prod_{i\in I} e_{n_i
\mu_i}$ insures that $x(e_{\mu})\neq 0$. But, as $\mu$ is an
interior point of $\Face$, there is also a relation $L \mu=\sum_{k,
\mu_k\in \Face} n_k\mu_k$ where $L$ and $n_k$ are strictly positive
integers. The relation $x(e_{\mu})^L=\prod_{k, \mu_k\in
\Face}x(e_{\mu_k})^{n_k}$ implies that $x_k\neq 0$ for all $k$ such
that $\mu_k\in \Face$. Thus $C(I)=\Face$.

\end{proof}

Let $\Face$ be a face of the cone $C$. Let $I(\Face)$ be the subset
of $\{1,2,\ldots,n\}$ such that $\mu_i\in \Face$. Let
$\ggy_{\Face}=\oplus_{i\in I(\Face)} \R \omega_i$ and
$\gr_{\Face}=\gr\cap \ggy_{\Face}$.

Assume $\gr_{\Face}$ is not equal to $\{0\}$. This means that
there exists a non trivial relation between the elements $\mu_i$
of $\tB$ belonging to the face $\Face$.
 For each $i\in
I(\Face)$, we still denote by $\nu_i$ the restriction to
$\gr_{\Face}\subset \gr\subset \ggy$ of the coordinate $\omega^i$
and we denote by $U_{\Face}:=\{u\in \gr_{\Face,\C}, \prod_{i\in
I(\Face)}\nu_i(u)\neq 0\}.$ For $m=\sum_{i\in I(\Face)} m_i\omega_i$
in $\ggy_{\Face,\Z}$, the restriction of the function
$p_m(z)=\prod_{i\in I(\Face)} z_i^{m_i}$ to $\gr_{\Face,\C}$ is well
defined on the open set $U_{\Face}$ of $\gr_{\Face,\C}$.

\begin{lemma}
Let $\Face$ be a face of $C$. If $\gr_{\Face}\neq \{0\}$, then
there exists a non-empty finite set $S_{\Face}\subset
\gr_{\Face,\Z}$ and constants $c_m \in \C^*$ such that the
relation $\sum_{m\in S_{\Face}} c_m p_m=0$ holds identically on
$U_{\Face}$.

\end{lemma}

\begin{proof}
Let us consider a basis
  $\{w_1, w_2,\ldots, w_{r_{\Face}}\}$ of $\gr_{\Face,\Z}\subset \ggy_\Z$.
  Here $r_{\Face}$ is the dimension of $\gr_{\Face}$.
  We denote simply by $p_j$ the rational function $p_{w_j}$ on
  $\C^n$. Each of the function $p_j$ is homogeneous of degree $0$.
  The rational map $p=(p_1,p_2,\ldots, p_{r_{\Face}})$ defines by
  restriction a map from $U_{\Face}$ to $(\C^*)^{r_{\Face}}$ which is homogeneous
  of degree $0$. Thus its image is of dimension strictly less than
  $r_{\Face}$. This implies that there is a non zero polynomial $Q$ such
  that $Q(p_1(u),\ldots,p_{r_{\Face}}(u))$ is identically $0$ on $U_{\Face}$.
  Using the relation $p_{w}p_{w'}=p_{w+w'}$, for $w,w'\in \gr_\Z$, the relation
  $Q(p_1(u),\ldots,p_{r_{\Face}}(u))$ is equivalent to a relation of the
  form given in the lemma.
\end{proof}

\begin{definition}
  For each face $\Face$ of the cone $C$ satisfying $\gr_{\Face}\neq \{0\}$, we
  choose a non-empty finite set $S_{\Face}\subset \gr_{\Face,\Z}$ and non zero
  constants $c_m, m\in
  S_{\Face}$, such that the rational function $R_{\Face}=\sum_{m\in S_{\Face}} c_m p_m$
  is identically equal to $0$ on $U_{\Face}$. We define
$$R=\prod R_{\Face},$$
where the product runs over the faces $\Face$ of $C$ with
$\gr_{\Face}\neq \{0\}$.

\end{definition}

\begin{definition}\label{Zsigma}
Let $\sigma$ be a subset  of  $\{1,2,\ldots, n\}$. Denote by
$$Z_\sigma:=\{z=(z_1,z_2,\ldots, z_n)\in (\C^*)^n\,;\, z_i=1, i\in
\sigma\}.$$
\end{definition}

\begin{lemma}
Let $\sigma$ be a subset  of  $\{1,2,\ldots, n\}$ such that the
set $\{\mu_i, i\in \sigma\}$  is a set of linearly independent
vectors of $V$. Then the function $R$ does not vanish identically
on $Z_\sigma$.
\end{lemma}

\begin{proof}
The function $R$ is of the form $\sum_{m\in \gr_\Z} k_m p_m$,
where $k_m$ is non zero for a finite non-empty set of $m\in
\gr_\Z$. Let
$$Y_\sigma=\{y=(y_1,y_2,\ldots, y_n)\in \C^n \, ; \,y_i=0, i\in
\sigma\}.$$ Denote by $e^y$ the  element
$(e^{y_1},\ldots,e^{y_n})$ of $Z_\sigma$. For $m\in \gr_{\Z}$, the
function $p_m$ is written as the exponential function
$$p_m(e^y)=e^{\sum_{k\notin\sigma}^n y_k \nu_k(m)}.$$

By Gale duality, if the vectors $\mu_i$ for $i\in \sigma$ are
linearly independent, then the linear forms $\nu_k$ for $k\notin
\sigma$ span $\gr^*$. For $m_1\neq m_2\in \gr_\Z$, the two linear
functions $ y\mapsto \sum_{k\notin \sigma}y_k \nu_k(m_1)$ and
$\sum_{k\in \notin\sigma}y_k \nu_k(m_2)$ are not identically equal
on $Y_\sigma$. Otherwise, we would obtain $\sum_{k\notin \sigma}
y_k \nu_k(m_1-m_2)=0$. This implies $m_1=m_2$ as $\sum_{k \notin
\sigma} y_k \nu_k$ ranges through $\gr^*$, when $y$ ranges through
$Y_\sigma$. Thus the relation $\sum_{m\in \gr_\Z} k_m p_m(e^y)=0$
between different exponentials cannot hold identically on
$Y_\sigma$.
\end{proof}

In particular, the set of $z\in Z_\sigma$ where $R(z)\neq 0$ is an
open dense set of $Z_\sigma$. As $R$ is the product of the
functions $R_{\Face}$ over the faces $\Face$ such that
$\gr_{\Face}\neq \{0\}$, it follows that the condition $R(z)\neq
0$ implies that $R_{\Face}(z)\neq 0$ for any face $\Face$  of $C$
such that $\gr_{\Face}\neq \{0\}$.

Finally, we obtain a  condition which ensures that  the sections
$\{f_{z,a_j},j=0,\ldots,d\}$ do not have common zeroes on
$\Toric(C)$.

\begin{prop}\label{denseopen}
Let $\sigma$ be a subset  of  $\{1,2,\ldots, n\}$ such that the
set $\{\mu_i, i\in \sigma\}$ is a set of linearly independent
vectors. Then for $z$ varying in the open dense set $Z$  of
$Z_\sigma$ where $R(z)\neq 0$, the $d+1$ sections $
f_{z,a_j},j=0,\ldots,d$ do not have a common zero in $\Toric(C)$.
\end{prop}

\begin{proof}
Consider a common zero $x$ of the functions $f_{z,a_j}$ on $\Aff(C)$.
We need to prove that $x_i=0$ for all $i\in \{1,2,\ldots,n\}$. If
not, there is a face $\Face$ of $C$, with $x_i\neq 0$ if and only
$\mu_i\in \Face$. The relation $\sum_{i=1}^n z_i x_i  \mu_i=0$
implies that $\gr_{\Face}\neq \{0\}$ for this face $\Face$. The
element $\sum_{i=1}^n z_i x_i \omega_i$  is in  $U_{\Face}$.

 Consider the rational function
$R_{\Face}=\sum_{m\in S_{\Face}} c_m p_m$ with $S_{\Face}\subset
\gr_{\Face,\Z}$. For $m\in \gr_{\Face,\Z}$, the value
$p_m(x_1,x_2, \ldots, x_n)$ is well defined and equal to $1$. Thus
the value of the function $R_{\Face}=\sum_{m\in S_{\Face}} c_m
p_m$ on the element $\sum_{i=1}^n z_i x_i \omega_i$ of $U_{\Face}$
is simply equal to
 $R_{\Face}(z)=\sum_{m\in S_F}c_m p_m(z)\neq 0$ by our assumption on $z$.
  This contradicts the fact that
$R_{\Face}$ vanishes identically on $U_{\Face}$.

\end{proof}

Thus, when $z$ is generic in $Z_\sigma$, the $d+1$ sections $
f_{z,a_j},j=0,\ldots,d$ satisfy the condition \eqref{cond11} which
ensures the existence of the toric residue.

Recall that we have chosen a  vector $\gamm\in V_\Z$ in the
interior of the cone $C$ and that  $\gwt \subset V^*$ is the
orthogonal of $\gamm$. The space $\gwt$ is a rational subspace of
dimension $d$ of $V^*$.
\begin{prop}\label{possible}
Let $[a_1,a_2,\ldots,a_d]$ be a basis of $\gt^*$. Let $\sigma$ be
a subset of $\{1,2,\ldots, n\}$ such that the set $\{\mu_i, i\in
\sigma\}$ is a set of linearly independent vectors. Then for $z$
varying in the open dense subset $Z$ of $Z_\sigma$ where $R(z)\neq
0$, the common zeros $x$ of the sections $f_{z,a_j}$ for
$j=1,\ldots,d$ lie in the torus $\Torus(H^*)$.
\end{prop}

\begin{proof}

As before, let $x\in \Aff(C)\setminus\{0\}$. If $x$ is a zero of the
functions $f_{z,a}$ for all $a\in \gwt$,  we want to prove that
$x_i\neq 0$ for all $i$.

From the equations $f_{z,a}(x)=0$, we obtain that
$$\sum_{i=1}^n z_i \langle a,\mu_i\rangle  x_i=0$$
for all $a\in V^*$ orthogonal to $\gamm$. Thus there exists $t\in
\C$ such that
$$\sum_{i=1}^n z_i x_i\mu_i=t \gamm.$$

If $t=0$, then the point $x$ is a zero of all sections $f_{z,a}$
for all $a\in V^*$ and the condition $R(z)\neq 0$ implies $x=0$.
So we may assume $t\neq 0$, and $\gamm$ belongs to the space
generated by the $\mu_i$ with $x_i\neq 0$. In particular there is
a relation $N\gamm=\sum_{i, x_i\neq 0} n_i\mu_i$ with $N>0$ and
$n_i\in \Z$. This implies that $x(e_{\gamm})\neq 0$.

 We have seen that
there exists a face $\Face$ of $C$ such that $x_i\neq 0$ if an only
$\mu_i\in \Face$. As $\gamm$ is an interior point of $C$, this face
is necessarily the cone $C$ itself. This concludes the proof of
the proposition.
\end{proof}

\end{document}